# Studies in Cryptological Combinatorics

by

Marc Zucker

This manuscript has been read and accepted for the Graduate Faculty in Mathematics in satisfaction of the dissertation requirement for the degree of Doctor of Philosophy.

__________________ _____________________________________

Date Chair of Examining Committee

__________________ _____________________________________

Date Executive Officer

Michael Anshel
Burton Randol
Al Vasquez

Supervisory Committee

THE CITY UNIVERSITY OF NEW YORK

# Abstract

## Studies in Cryptological Combinatorics

by

Marc Zucker

Advisor: Professor Michael Anshel

The key-agreement problem (finding a private key to use for secret messages, otherwise referred to as the public-key distribution problem), was introduced by Diffie and Hellman in 1976. An approach to structuring key-agreement protocols via the use of one-way associative functions was proposed in 1993 by Rabi and Sherman. We propose here a provably strong associative one-way function based upon knot composition (answering an open problem proposed by Rabi and Sherman whether any such associative one-way functions exist).

We also introduce and solve a game, exploring its relation to problems in graph and braid theory and develop a new technique for computing whether a graph is *n*-colorable. En route we look at estimator and prediction problems raised in Classical Probability Theory using Urn problems.

## Acknowledgements:

I would like to give special thanks to my advisor, Michael Anshel, for his guidance and support.  I would also like to thank Joseph Malkevitch for his helpful conversations and suggestions.

# Table of Contents

## List of Figures

# 1 Outline and Summary

The key agreement problem, that of finding a private key to use for secret messages (otherwise referred to as the public-key distribution problem), was introduced by Diffie and Hellman in 1976 [18]. The underlying security of these schemes is usually based upon some computationally 'hard' problem. A problem is considered hard if it is not computable in polynomial time. (The Diffie-Hellman key agreement scheme was based upon the then believed idea that discrete logarithms over $GF(p)$ is difficult to compute.) Several other approaches to this problem have come out since then. One of them, an approach to structuring key-agreement protocols, is with the use of one-way associative functions, introduced in 1993 by Rabi and Sherman [51]. Public-key cryptography and key-agreement protocols using algebraic techniques, has had much recent development as seen from recent papers by Iris Anshel, Michael Anshel, and Dorian Goldfeld (among others) [4, 5, 6]. We introduce an associative one-way function based upon knot composition that under certain given intractability conditions we consider strong (answering a question proposed in [51] of whether any such associative one-way functions exist). An associative one-way function is considered strong if and only if given either the first or the second argument inverting the function is not computable in polynomial time. This we use to show the primary idea presented in this paper; that is, of the existence (and an example) of a proposably strong secret key-agreement protocol.

One of the classical techniques in code breaking has been that of the statistical attack. Though of little use in modern cryptography, the underlying role that probability has played in cryptographic systems, whether in the making or breaking of said systems, is still present. Classical probability started as a means of modeling real life situations. Initially this was in the realm of games of chance, but it developed later in a number of different directions. Polya, for instance, modeled the spread of viruses throughout a civilization with his famous urn model. The use of urns as a means of modeling was used as early as Laplace [40], using this model to derive a law of succession. That is, given an urn with an unknown distribution of white and black marbles Laplace derived a law that would tell us what the probability of pulling a white marble from the urn would be, given that every other time a marble had been pulled out of the urn (and then returned) it had been white. A great number of debates have ensued about the validity of Laplace's law due to the assumptions that must be made about the prior probabilities. The notion of a successor, though, is one that is not only fascinating in its own right, but one that has natural interest for the cryptographer. And it was, in fact, these same ideas that were at play behind an estimator developed by Alan Turing and I. J. Good used to try to break the Enigma cipher during World War II. While there are only speculations about the underlying ideas behind this estimator, it has aroused interest into the nature and function of estimators. While the Good-Turing estimator was used for cryptanalysis, the theory and application of estimators has found important uses in data communications, image reconstruction, and language recognition. With this information and the underlying theory of Laplace as a start, we look at how Laplace's law is affected by differing priors as well as ongoing conditions and what happens as the underlying distribution becomes unbounded.

Studies in mathematical games have become a tool leading to a wealth of new information regarding different mathematical structures. As already mentioned, ever since Pascal and Fermat created the classical theories of Probability games have been studied as a science. However games qua games have only recently been studied and used as models (along with, if not developed from, probability) for differing situations throughout many different fields. Game Theory, introduced with the publication of *Theory of Games and Economic Behavior* by Von Neumann and Morgenstern [48], not only showed that the study of games (and strategy) can be viewed as a branch of mathematics, but that this field can be used to model realistic situations with great success. Ever since the publication of that seminal text, Game Theory has been used to do just that in areas such as economics, social behavior, and warfare, among others. Sadun, Villegas, and Voloch [54] introduced the one-person game of Blet, and with it came its ability to model the 3-Braid. The relation between braids and games is one that has not been sufficiently investigated and has therefore motivated us to investigate games as a means of studying various algebraic and geometric structures.

# 2 Knot Composition and Secret Key-Agreements

## 2.1. Introduction

The use of one-way associative functions as a paradigm for secret key-agreement and digital signatures was first introduced by Rabi and Sherman in 1993 [51]. It was left as open problems to exhibit a plausible strong associative one-way function and a provable strong one-way function. We answer these questions via the use knot composition [57], which is considered plausibly strong.

Knots and braids have been mathematical structures of interest for a long time. The classification and determination of the equivalence of knots, are two classic problems in the theory knots. While the study of knots as a mathematical subject is not that old, it predated the study of braids as such. The 19th Century saw the study and development of knots with Lord Kelvin and P. G. Tait (among others), while the 20th Century saw further advancements with the achievements of K. Reidemeister and J. Alexander. Braids, however, though looked at as early as Gauss [23], first received its first real formulation as a mathematical object with the papers of Emil Artin [7, 8] early in the 20th Century.

Algebraic, and more specifically group theoretic, methods as used in cryptography, and more importantly public-key cryptography, have been developing at ever more rapid rates. See Anshel, Anshel, Goldfeld, and Fisher ([3, 4, 5, 6]) for a development of such ideas including the use of key-agreement protocols within an algebraic setting. The use of knots and braids in cryptographic schemes has therefore been receiving more and more attention as more people study problems in both of these areas. One outcome of the increased interest in this area is the recognition of some problems as (computationally) hard (such as, for example, the Conjugacy Search Problem in the braid group [32]). For this reason cryptologists have naturally turned their attention toward the study of knots and braids.

Problems, inherent and accidental, to knots and braids have therefore become a feeding ground for new ideas and techniques in cryptography. The braid conjugacy problem is one of the more known problems with applications in cryptography. (For a recent paper on specific case of positive braid permutations and conjugacy see Morton and Hadji [44].) But it is not simply cryptologists, but mathematicians and scientists in general, who have turned their attention to these areas, for its applications and importance reach all areas; from physics and the motion of celestial bodies to biology and the knotting's of DNA.

Since it is possible for us to represent knots as the closure of braids, and braids are more easily represented, they have become a natural field of study for knot theorists. We therefore focus our attention as well on knot composition and secret key-agreement protocols as represented by braids.

## 2.2. A Survey Of Cryptological Security

The security of a system is the integral component in any cryptological scheme. What secure means, however, has changed over time as different and newer techniques of breaking codes have been developed. For Caesar the security was in a permutation. Statistical analysis has shown the lack of security involved in that system. Turing and Good used similar techniques to try to break the German's Enigma code.

What is meant by the security of a system now comes under a couple of different definitions based upon the model of the system. Most commonly it is a level of the mathematical complexity of a system. Systems can be such as those in which the scheme is known only to the individual parties, one in which insecure lines of communication are used primarily to establish a secret key agreement between the parties, or one known as a public-key system.

In 1974 G. Purdy [50] introduced the first detailed one-way function. A one-way function is one in which encryption is considered simple while reversing the process (decryption) is considered hard. 1976 found the introduction of the Diffie-Hellman protocol for establishing public keys [18]. Public-key cryptosystems introduced a new paradigm for cryptology – one in which a trapdoor function exists. The complexity of the system thus rested upon the difficulty of breaking this trapdoor. In 1977 Rivest, Shamir, and Adleman developed what is known as the RSA cryptosystem, perhaps the most popular of the public-key cryptosystem [53].

But as new systems get developed newer techniques for breaking them are developed. As well, the security of any system is based upon certain defined and understood conditions; this is especially the case in public-key cryptosystems in which the security has commonly been based upon the complexity of a given problem. (As, for example, the discrete logarithm within the Diffie-Hellman scheme.) RSA, based upon the difficulty of integer factorization has been under attack ever since the system was first introduced. Boneh and Venkatesan [14] have suggested that attacks other than that of integer factorization might be possible upon RSA. See [12, 35] for an overview of various different attacks upon RSA as well upon other systems that have been developed. (Systems based upon the Knapsack problem, as well as Elliptic Curves and other recent approaches are presented in this paper.) Also see [25] for a general overview of cryptographic security.

Recently it has become increasingly more apparent that other means of attack are possible. These types of attacks, known as 'side-channel' attacks take advantage of seemingly extraneous information stemming from the operating environment and specific properties of the implementation of the system such as the execution time [36], the amount power consumed [37], error messages [11, 42], induced errors [13], and electromagnetic radiation [1.

Recently braid groups have been used as the basis for cryptological schemes [17]. Most prominently has been the use of the conjugacy problem in braid groups to build a system. For an overview of the use of braid groups in cryptography and in specific the use of the conjugacy problem see [4, 32, 33, 34].

With the introduction of the notion of Associative One-Way Functions by Rabi and Sherman in 1993 [51] new paradigms in protocols were being established. In 1997 [52]

Rabi and Sherman showed that such functions exists if and only if $P \neq NP$. See [26] for a discussion of whether it is possible at all to base the security of a system on $P \neq NP$. Homan, in 2000, gave an upper bound for the ambiguity of these functions [30]. Further developments in security protocols have been established with the introduction of Zero-Knowledge Proofs [24] and more recently with the Random Oracle Model [16]. Even more recently the advent of Quantum Cryptology has tested the limits and security of cryptological systems. See [41, 56].

## 2.3. Associative One-Way Functions

Rabi and Sherman [51] proposed in 1993 the idea of associative one-way functions as a paradigm for secret-key agreement protocols. Throughout this paper we assume that there are two parties, Alice and Bob, which wish to send information to each other that is to be kept secret. We assume, as well, that a third party, whom we call Eve, intercepts their transmissions. The ability to set up a secret-key agreement between the parties is therefore crucial to being able to communicate without the fear of the communication falling into the hands of those they do not wish to obtain such information. Associative one-way functions are set up as a paradigm for creating secret-keys by which a cryptographic scheme can be created with assured security (within reasonable limits). We present here a basic outline of associative one-way functions. (See [51] for a more thorough presentation.)

**Definition 2.3.1** A binary function, $\circ : S \times S \to S$, is said to be *one-way* if and only if $\circ$ is honest and computable in polynomial time, whereas inverting $\circ$ is not computable in polynomial time.

Rabi and Sherman [51], deal exclusively with $S$ as binary functions on the infinite message space $S = \{0,1\}^*$ of all finite binary strings. (A function, $\circ$, is said to be honest if and only if there exists a polynomial $p$ such that, for every $z \in \text{image}(\circ)$, there exists $x, y \in S$ such that $x \circ y = z$ and $|x| + |y| \leq p(|z|)$. The necessity for the honesty requirement arises so that the complexity (difficulty) of inverting an associative one-way function is not based solely upon the input strings being significantly larger than the output string.)

**Definition 2.3.2** Any binary one-way function, $\circ : S \times S \to S$, is said to be *strong* if and only if given its first argument inverting the function is not computable in polynomial time. And similarly if given its second argument.

**Definition 2.3.3** A (binary) function is said to be *associative* if and only if $x \circ (y \circ z) = (x \circ y) \circ z$ holds true for all $x, y, z \in S$.

**Definition 2.3.4** Any binary function, $\circ : S \times S \to S$, is said to be a *(strong) Associative One-Way Function* if and only if $\circ$ is both associative and (strongly) one-way.

Integer multiplication is an example of an associative one-way function, since the splitting of an integer is considered hard (non-polynomial time). However, since division is simple given one of the factors, it is not a strong associative one-way function. Similarly the Diffie-Hellman key exchange algorithm for creating secret-keys is not associative, and therefore not an associative one-way function. As well, the hardness of both of these problems is questionable.

Rabi and Sherman leave it as an open problem to find a *provably strong* associative one-way function. Where by *provably strong* was meant: an implementation that under suitable intractability assumptions is shown to satisfy the definition of a strong associative one-way function. In this paper we shall use the phrase *proposably strong* instead of their *provably strong*. We propose here a proposably strong associative one-way function as a protocol for secret-key exchanges. This protocol is based upon the assumed intractability (hardness) of knot factoring.

## 2.4. Definition of the Braid Group

We present in this section an overview of what a knot is, what a braid is, and what the braid group is.

**Definition 2.4.1** A Knot is a one-dimensional closed non-intersecting curve in three-dimensions. (We can also say that it is a (smooth) embedding of $S^1 \to S^3$ or $R^3$.)

The simplest knot is the un-knot, that is, the circle in three dimensions.

Classically, knots are studied by studying their projection in two-dimensions. Two projections can be shown to represent the same knot if there is a set of moves, called Reidemeister moves, that changes one into the other.

**Definition 2.4.2** A Braid is a set of one-dimensional, non-intersecting, curves in three dimensions all of which descend strictly monotonically from one given horizontal line to some other given horizontal line below.

In general, we will leave the horizontal lines out of our drawings. (Figure 2.1)

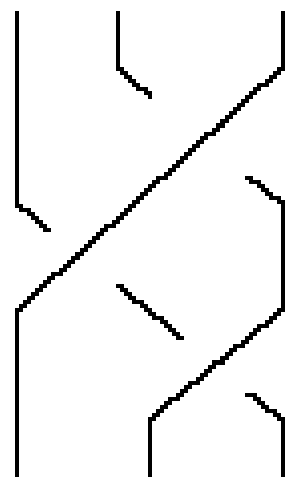

Figure 2.1

Two braids are considered equivalent if the strands of one braid can be manipulated, without any of them intersecting any other, so that it is made to look like the other braid. (Figure 2.2)

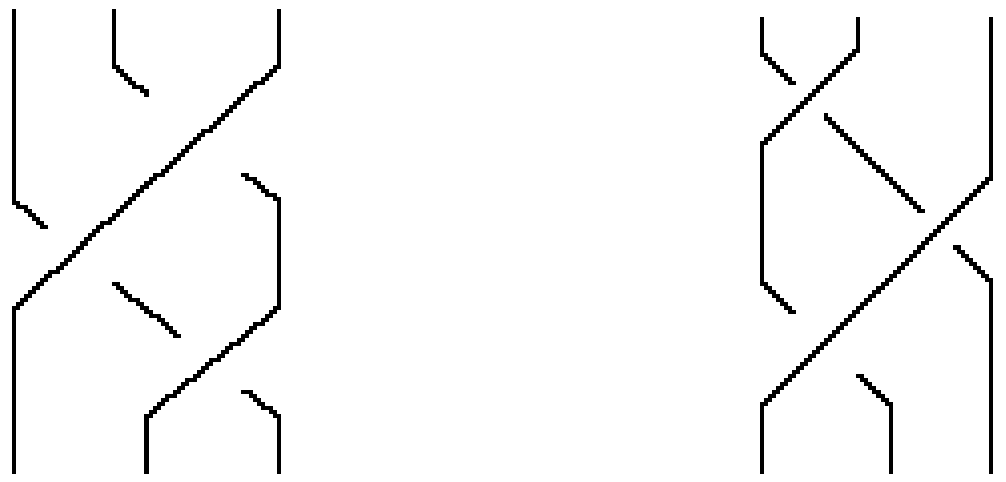

Figure 2.2 Two equivalent braids on three strands.

The class of all braids with a given number of strands forms a group under the operation of braid addition. The addition of braids is formed by joining the top of the strands of the second braid to the bottom of the strands of the first braid (and removing the horizontal lines/bars if used). For example,

composing 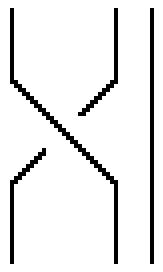 with 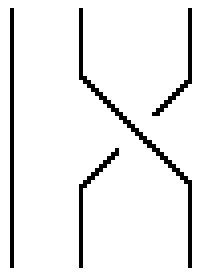 yields. 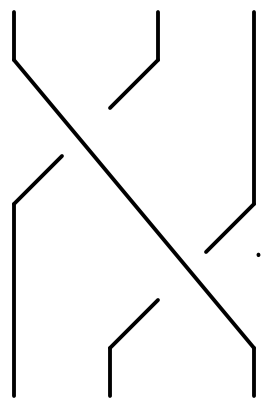

The inverse of a braid, found by taking the braid's mirror image (from the bottom up), is seen in figure 2.3 and can be easily checked. Lastly, the group's associativity arises from its very construction. The necessary components of a group have therefore all been accounted for. In fact the different classes of braids, i.e. the classes of different numbers of strands, form equivalence classes of braid groups.

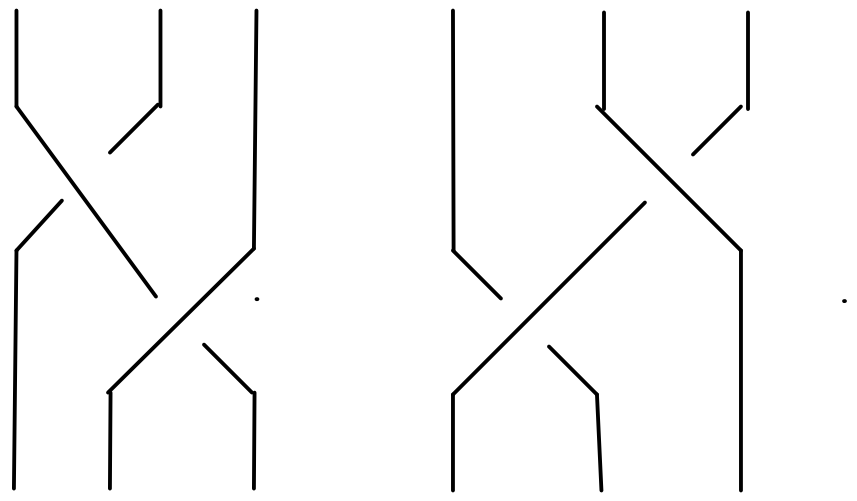

Figure 2.3 A braid and its inverse.

## 2.5. The Artin Presentation of the Braid Group

Emil Artin introduced in the early 20th century not only a group structure in braids, but also a specific presentation of one in terms of generators and relators [7, 8]. This presentation was as follows: take horizontal lines (or planes if we are looking at the three dimensional version) and divide a braid in such a way so that between any two parallel (horizontal) lines there is at most only two strands which cross at most once. If between any two lines there are no crossings, then we can remove this section without any loss of the braid structure. Each individual section of this braid can now itself be viewed as a braid and the braid as a whole as the sum of all these individual braids.

Let the strands of the braid be labeled from left to right, and let $\sigma_i$ indicate that braid in which only the $i^{th}$ and $(i+1)^{st}$ strands cross (once) with the $(i+1)^{st}$ strand going over the $i^{th}$ strand. (Figure 2.4.)

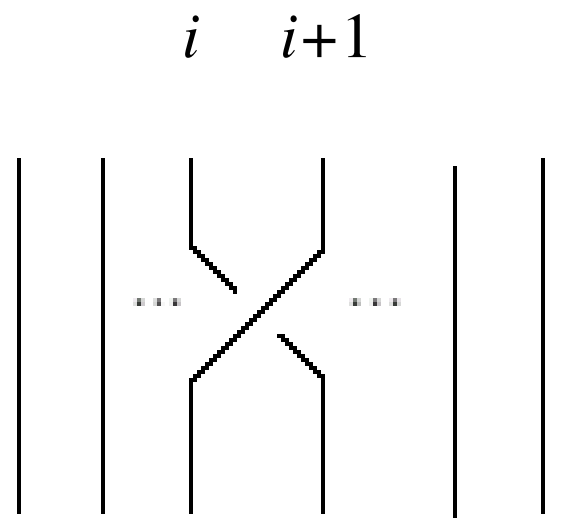

Figure 2.4

The inverse of $\sigma_i$ is indicated by $\sigma_i^{-1}$, the braid in which the $i^{th}$ and $(i+1)^{st}$ strands cross, but with the $i^{th}$ strand going over the $(i+1)^{st}$ strand instead.

The presentation (i.e. a way of specifying a group completely via generators and relators), of $B_n$, the Braid Group of index *n* (that is the braid group on *n* strands), given by Artin was defined with the following generators and relators.

**Generators:** $\sigma_1, \sigma_2, \ldots, \sigma_{n-1}$

**Relators:** $\sigma_i \sigma_j = \sigma_j \sigma_i \quad |i - j| > 1$

$\sigma_i \sigma_{i+1} \sigma_i = \sigma_{i+1} \sigma_i \sigma_{i+1}$

(Figure 2.2 shows an example of the second relation.)

## 2.6. Knot Composition

### 2.6.1. Connected Sums

What we would like to achieve now is a notion of adding knots. (We will deal only with oriented knots.) One way is to cut each knot and attach one loose end of each together, and then to do the same with the other loose end (making sure to match up the orientations).

For example, to add

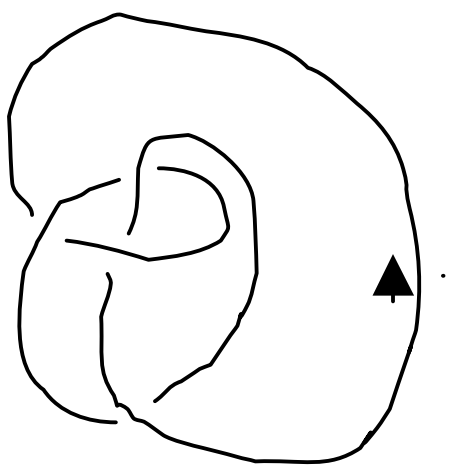 and 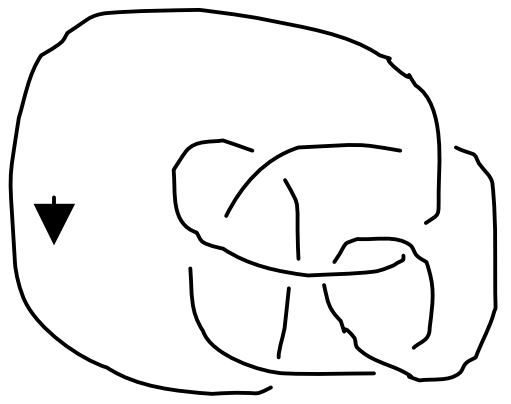 we cut each knot as such

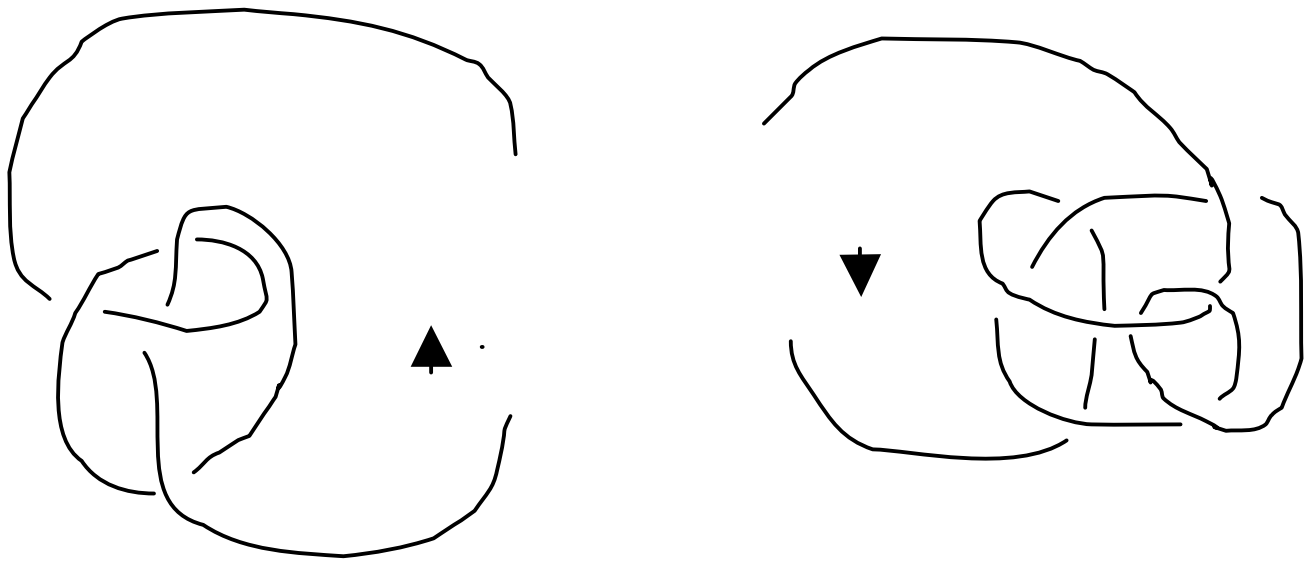

And then attach the loose ends, resulting in:

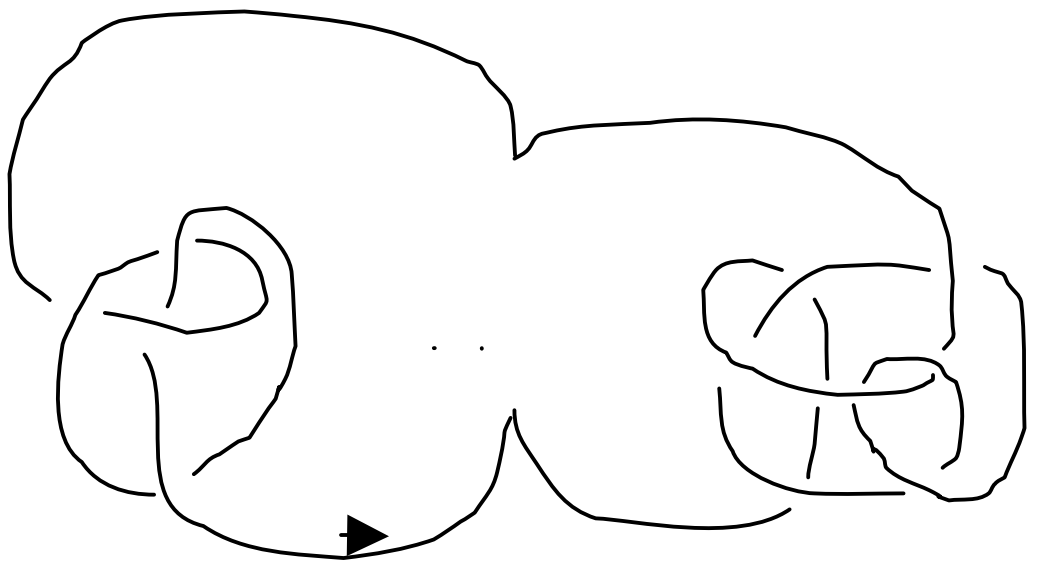

Let us now formalize this notion of knot addition (called knot composition), and its opposite, knot factoring, following Sullivan [57].

**Definition 2.6.1** Let $K_1$ and $K_2$ be knots in distinct 3-spheres $S_1^3$ and $S_2^3$, respectively. Pick two points $a_1 \in K_1$ and $a_2 \in K_2$, and choose two small balls $B_1$ and $B_2$, centered at $a_1$ and $a_2$, respectively, such that $B_i \cap K_i$ can be deformed to an axis of $B_i$, for $i$ = 1,2. Form a union $S_1^3 \setminus \overset{o}{B}_1 \cup S_2^3 \setminus \overset{o}{B}_2$, using a gluing homomorphism that matches $K_1 \cap \partial B_1$ to $K_2 \cap \partial B_2$ with the exiting endpoints going to the entering endpoints. Thus we have a new 3-sphere containing a new knot called the *connected sum* of $K_1$ and $K_2$, which we denote by $K_1\#K_2$.

The reverse of definition 2.6.1 is called *knot factoring*. It can be shown that this definition for the connected sum of two knots is independent of the choice of points made. It follows from here that connected sums are commutative, that is $K_1\#K_2 = K_2\#K_1$. This can be seen descriptively by allowing one knot to become extremely small and allow it to '*slide*' along the knot until it is in the position desired, at which point it can be enlarged to its original size. (See [57] for a full discussion of knot factoring.) Associativity, $K_1\#(K_2\#K_3) = (K_1\#K_2)\#K_3$, is also easy to show.

## 2.6.2. Braid Closures as Knots

Given an arbitrary braid we can form what we call the closure of the braid by adjoining the $i^{th}$ loose strand at the top of the braid with the $i^{th}$ loose strand at the bottom of the braid. I.e. As if the bottom of the braid was added (glued) to the top of the braid. (Figure 2.5.) Depending upon the given braid its closure might be a knot, the unkot, or a set of links. It can be shown that given any knot it is possible to find a braid such that the closure of this braid is the given knot (see [47, 43, 15, 46, 31]). We refer to any braid whose closure is a given knot as a braid representation of that knot. Figure 2.6 shows a knot and a braid representation of that knot.

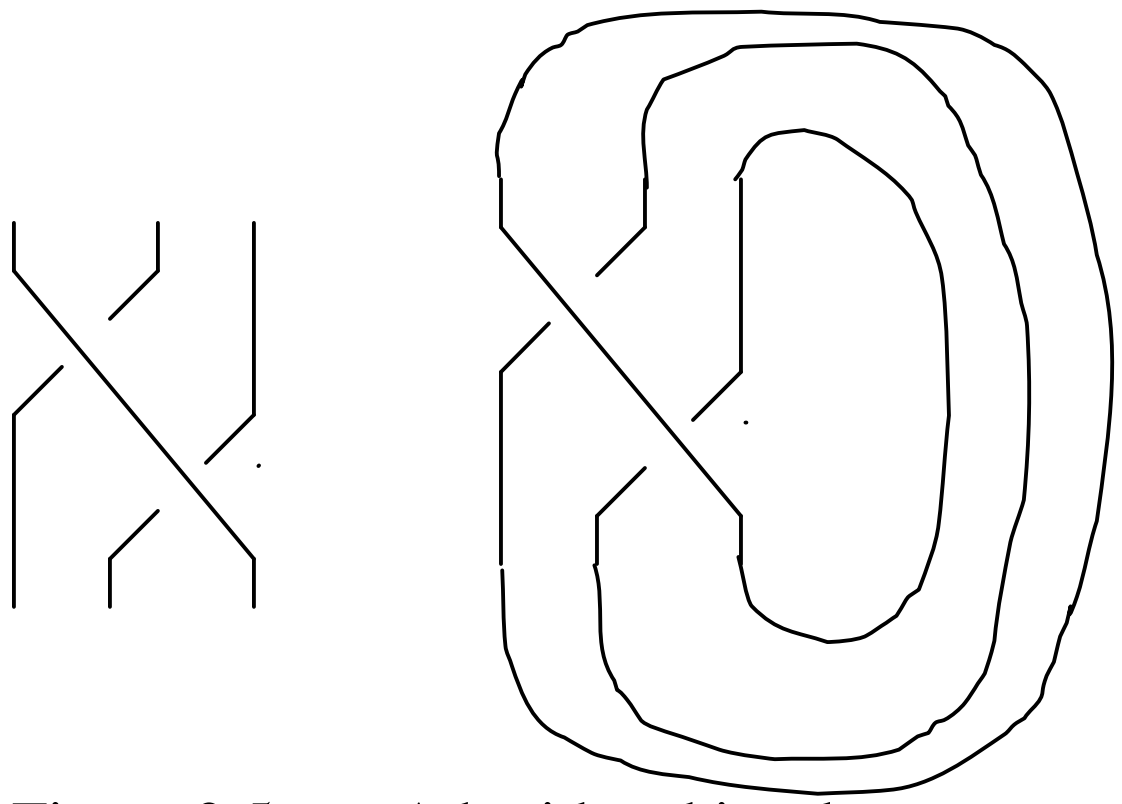

Figure 2.5 A braid and its closure.

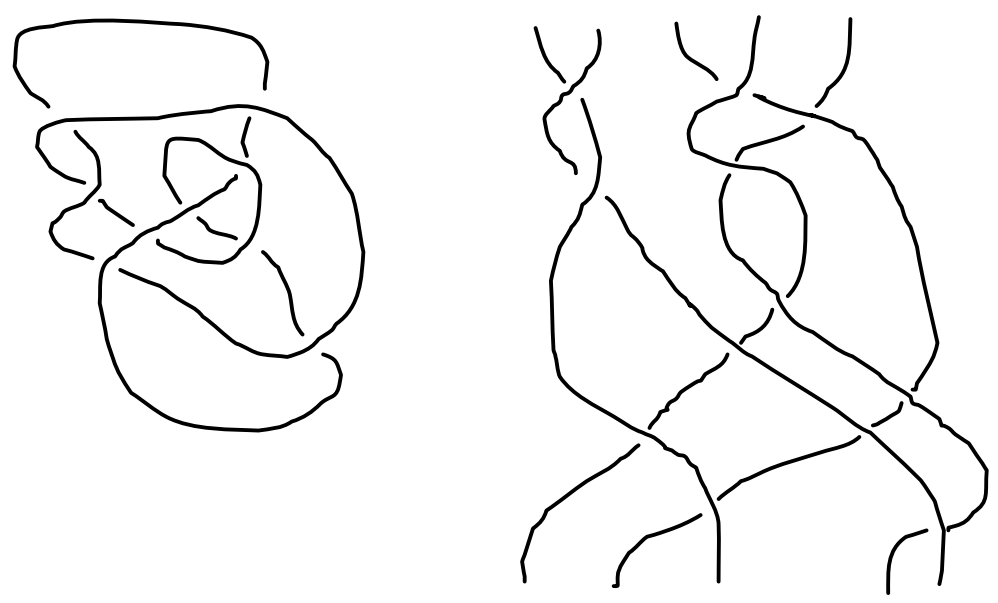

Figure 2.6 A knot and a braid whose closure is that knot.

Let us present, in outline form, how a knot might be transformed into a braid (this basically follows Murasugi and Kurpita [47]).

On an oriented knot let us mark off points such that the segments between alternating points alternate go from an overcrossing to an undercrossing. (The segments are allowed to have no crossings.) We now have *2m* points marked off. Next, let us place circles around every other point that we just marked off; in specific those that start an overcrossing segment. Let us now join these circles by bands such that: (i) each band meets each circle, if at all, at only a single edge of the band, (ii) a band does not contain any crossing points of the knot, and (iii) the band does not contain any points that we labeled other than those in the circles.

We can now imagine that this new loop (created by the connected bands, for they all join circles) that goes through this knot is straightened out (that is, flattened out so as to lie in a plane). It should be noticed not only that through this loop passes the strands of the knot, but that these strands all pass through in the same direction.

Let us now line up these strands, cut the parts of the strands where they go through the loop, and extend these edges upward and downward forming the top and the bottom of a braid. We have thus formed a braid whose closure is our original knot.

**Definition 2.6.2** Two braids are called *knot-equivalent* if their closures represent the same knot.

In the 1930's Markov stated that the closure of two braids represent the same knots if one can be changed into the other by two types of moves. These moves are now known as Markov moves. Thus, equivalently, two braids are said to be *knot-equivalent* if one can be transformed into the other via a series of Markov moves. Where the Markov moves are as follows:

- A Markov move of the first type is where a braid is replaced by its conjugate. That is where a braid $\beta$ is replaced by $\gamma\beta\gamma^{-1}$, where $\gamma$ is an arbitrary braid.
- A Markov move of the second type is where an $n$-braid, $\beta$, is replaced by an $(n+1)$-braid, $\beta\sigma_n$ or $\beta\sigma_n^{-1}$.

(Markov's proof was incomplete; a complete proof first appeared much later. See [47, 43] for a proof Markov's Theorem.)

# 2.7. Knot Composition via Braid Representations

We now present an algorithm for composing two knots via their braid representations.

Let $m$ and $n$ be the number of strands of $B_1$ and $B_2$, respectively, and let $<B_1>$ and $<B_2>$ be the words of those two braids (that is, the braid presentations of these two braids in terms of Artin's presentation). We compose the two braids by adjoining one strand from the end of the first braid with one strand from the beginning of the second braid while leaving the rest unchanged (figure 2.7). The composition of the two knots in braid representation is thus effectuated as follows:

$(B_1 \# B_2)=$

$$<B_1> \sigma_n^{-1}\sigma_{n-1}^{-1}\cdots\sigma_2^{-1}\sigma_{n+1}^{-1}\sigma_n^{-1}\cdots\sigma_3^{-1}\cdots\sigma_{m+n-1}^{-1}\sigma_{m+n-2}^{-1}\cdots\sigma_n^{-1} <B_2> \sigma_n^{1}\sigma_{n-1}^{1}\cdots\sigma_2^{1}\sigma_{n+1}^{1}\sigma_n^{1}\cdots\sigma_3^{1}\cdots\sigma_{m+n-1}^{1}\sigma_{m+n-2}^{1}\cdots\sigma_n^{1}$$

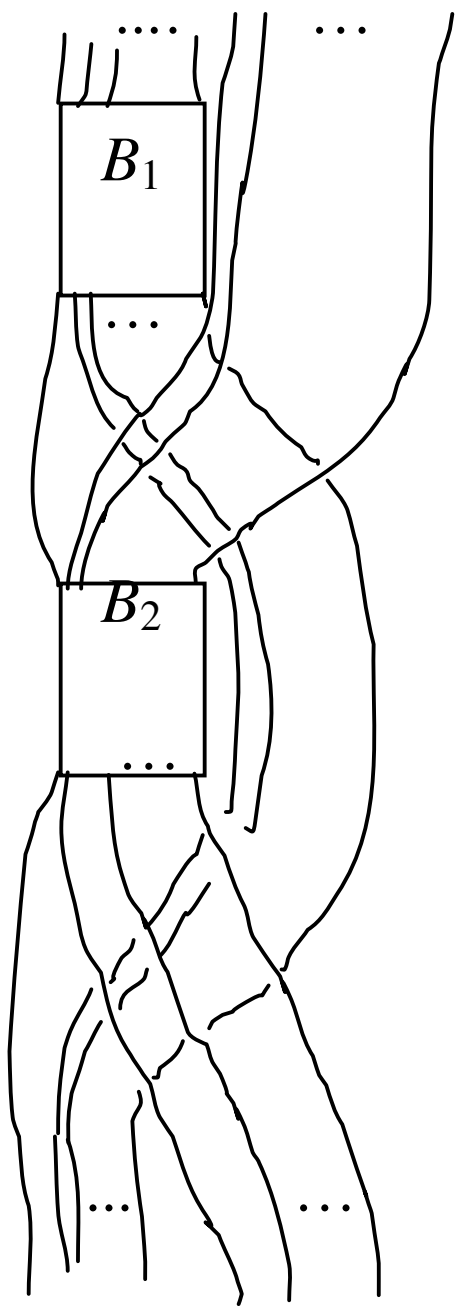

Figure 2.7 The composition of two braids

Similarly, if $B_1 = A_1 + A_2$ (where the addition is the standard addition of braids presented in 2.2), then we can also compose the braids as follows:

$(B_1\#B_2)=$

$$< A_1 > \sigma_n^{-1}\sigma_{n-1}^{-1}\cdots\sigma_2^{-1}\sigma_{n+1}^{-1}\sigma_n^{-1}\cdots\sigma_3^{-1}\cdots\sigma_{m+n-1}^{-1}\sigma_{m+n-2}^{-1}\cdots\sigma_n^{-1} < B_2 > \sigma_n^{1}\sigma_{n-1}^{1}\cdots\sigma_2^{1}\sigma_{n+1}^{1}\sigma_n^{1}\cdots\sigma_3^{1}\cdots\sigma_{m+n-1}^{1}\sigma_{m+n-2}^{1}\cdots\sigma_n^{1} < A_2 >$$

(This, of course, adds a little more complexity to the knot factoring that must be achieved by Eve in order to find the secret key.)

## 2.8. Secret Key-Agreement via Knot Composition

An associative one-way function can now be established as follows. Let Alice take two braids, $B_1$ and $B_2$, compose them to get $B_1\#B_2$, and then send $B_1\#B_2$, and $B_2$ to Bob. Bob then takes his own braid, $B_3$, composes it with $B_2$, getting $B_2\#B_3$, and sends it to Alice. (Remembering that knot composition is commutative and therefore the order of the composition is unimportant.) Alice then computes $B_1\#(B_2\#B_3)$ while Bob computes $(B_1\#B_2)\#B_3$. The two knots corresponding to the closure of these braids being equal, Alice and Bob can then create a secret key from computable knot invariants. Thus, for example, it is possible to compute the Jones polynomial (see [47]) and use some agreed upon coefficient as the secret-key.

Since knot factoring is considered intractable in polynomial time (i.e. hard), this function satisfies the requirements for a proposably strong associative one-way function. However, it will be necessary that certain operations be done to the braid word of the composition so that it is not obvious what the factors are. These can be accomplished via four operations/relations which leave the knot represented by the closure of the braid unchanged. These are the two braid identities: $\sigma_i\sigma_j = \sigma_j\sigma_i$ $|i - j| > 1$ and $\sigma_i\sigma_{i+1}\sigma_i = \sigma_{i+1}\sigma_i\sigma_{i+1}$, and the two Markov moves.

Given two braids it is known to be hard to determine whether or not they are the same. Thus the inclusion of these extra moves establishes the complexity needed to make knot factoring hard and therefore our secret-key agreement protocol as proposably strong.

A key-agreement protocol (along the lines of one presented in [51]), can be effectuated with any amount of parties by having one party select a braid, *A*, at random, and sends it to all other parties. Each party then selects a braid, $B_i$, at random, computes $A\#B_i$, and sends it to all the other parties. Each party then computes $B_i\#(A\#B_1)\#(A\#B_2)\#\ldots\#(A\#B_n)$ $= B_1\#B_2\#\ldots\#B_n\#A^{n-1}$, where $A^{n-1}$ means $(A\#A\#\ldots\#A)$ with # applied $n$-2 times. Due to the commutativity and associativity of knot composition, the closure of these braids will all be equivalent and a key can be established by some agreed upon coefficient of a knot invariant (or any other number or method of a similar sort). Or simply this coefficient can be used to confirm the identity of the two parties involved.

We have shown the following:

**Theorem 2.8.1** There exists a proposably strong associative one-way function. Furthermore, a secret key-agreement protocol can be established via the use of this associative one-way function.

## 2.9. A Simple Example

Let us assume that Alice and Bob have exchanged and composed their knots as we have described, and that their resulting knot is the following:

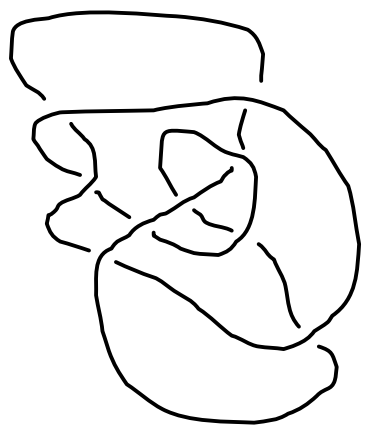

This knot, together with a braid whose closure was this knot, was given in figure 2.5.

We should note that the knot that Alice and Bob each compute will (most likely) have different representations. Our knot is just one possible representation that it might have. (It should also be noted that our example is an extreme simplification in so far as the fact that our knot is in fact prime, and therefore could not be the result of any composition. We are therefore simply using it to demonstrate the underlying technique.)

The Artin presentation for this braid is as follows: $\beta = \sigma_1\sigma_3\sigma_4^{-1}\sigma_3^{-1}\sigma_1\sigma_3^{-1}\sigma_2^{-1}\sigma_4^{-1}\sigma_1^{-1}\sigma_3^{-1}\sigma_4^{-1}\sigma_2^{-1}$.

Introduced by Vaughan Jones in 1984, the Jones Polynomial, which is invariant under different representations (and therefore of interest to us), is derived by the following formula:

$$\xi(\beta) = \frac{t^{\frac{1}{2}(\varepsilon(\beta)-n+1)} tr\left(\Phi_n(\beta)\mu^{\otimes n}\right)}{1+t},$$

where $\Phi_n\left(\sigma_i^{\varepsilon}\right) = \underbrace{I_2 \otimes \cdots \otimes I_2}_{i-1} \otimes R^{\varepsilon} \otimes \underbrace{I_2 \otimes \cdots \otimes I_2}_{n-i-1}$, $\mu = \begin{bmatrix} 1 & 0 \\ 0 & t \end{bmatrix}$, $\varepsilon(\beta) = \exp(\beta)$,

$R = \begin{bmatrix} 1 & 0 & 0 & 0 \\ 0 & 0 & -\sqrt{t} & 0 \\ 0 & -\sqrt{t} & 1-t & 0 \\ 0 & 0 & 0 & 1 \end{bmatrix}$, and $\otimes$ is a tensor product.

A *tensor product* of two matrices is defined as follows:

**Definition 2.8.1** Given matrices $\mathbf{A} = (a_{ij})$ and $\mathbf{B} = (a_{kl})$, respectively of size $p \times q$ and $r \times s$, then the tensor product of **A** and **B, denoted by** $\mathbf{A} \otimes \mathbf{B}$, is the $pr \times qs$ matrix defined as $\mathbf{A} \otimes \mathbf{B} = \begin{bmatrix} a_{11}B & a_{12}B & \cdots & a_{1q}B \\ a_{21}B & a_{22}B & \cdots & a_{2q}B \\ \vdots & \vdots & \ddots & \vdots \\ a_{p1}B & a_{p2}B & \cdots & a_{pq}B \end{bmatrix}$.

For example, if we were given the two matrices $\mathbf{A} = \begin{bmatrix} 0 & -3 \\ 2 & 1 \\ 1 & 1 \end{bmatrix}$ and $\mathbf{B} = \begin{bmatrix} 1 & -1 & 2 & 1 \\ 2 & 0 & 3 & 1 \end{bmatrix}$,

the tensor product of **A** and **B** would be:

$$\mathbf{A} \otimes \mathbf{B} = \begin{bmatrix} 0 \cdot B & -3 \cdot B \\ 2 \cdot B & 1 \cdot B \\ 1 \cdot B & 1 \cdot B \end{bmatrix} = \begin{bmatrix} 0 & 0 & 0 & 0 & -3 & 3 & -6 & -3 \\ 0 & 0 & 0 & 0 & -6 & 0 & -9 & -3 \\ 2 & -2 & 4 & 2 & 1 & -1 & 2 & 1 \\ 4 & 0 & 6 & 2 & 2 & 0 & 3 & 1 \\ 1 & -1 & 2 & 1 & 1 & -1 & 2 & 1 \\ 2 & 0 & 3 & 1 & 2 & 0 & 3 & 1 \end{bmatrix}.$$

Purely for the sake of simplification, so as to compute the Jones polynomial more easily (and so as to be able to represent it in a limited amount of space), we can apply the relations of the Artin presentation together with Markov moves to get $\beta = \sigma_1\sigma_3\sigma_4^{-1}\sigma_3^{-1}\sigma_1\sigma_3^{-1}\sigma_2^{-1}\sigma_4^{-1}\sigma_1^{-1}\sigma_3^{-1}\sigma_4^{-1}\sigma_2^{-1} = \sigma_2^{-1}\sigma_3^{-2}\sigma_2^{-2}\sigma_3^{-1}\sigma_1^{-1}\sigma_2^{-1}\sigma_1^{2}\sigma_3$. (Though a more reduced form might be possible, our objective here was simply to reduce it enough so as to make the presentation reasonable.)

Using the following matrices we compute $tr\left(\Phi_4(\beta)\mu^{\otimes 4}\right)$. (Since our braid has been reduced to a 4-braid (from a 5-braid as it was previously represented), we need only to compute this value for the 4-braid.)

$$\sigma_1 = \begin{bmatrix}
1 & 0 & 0 & 0 & 0 & 0 & 0 & 0 & 0 & 0 & 0 & 0 & 0 & 0 & 0 & 0 \\
0 & 1 & 0 & 0 & 0 & 0 & 0 & 0 & 0 & 0 & 0 & 0 & 0 & 0 & 0 & 0 \\
0 & 0 & 1 & 0 & 0 & 0 & 0 & 0 & 0 & 0 & 0 & 0 & 0 & 0 & 0 & 0 \\
0 & 0 & 0 & 1 & 0 & 0 & 0 & 0 & 0 & 0 & 0 & 0 & 0 & 0 & 0 & 0 \\
0 & 0 & 0 & 0 & 0 & 0 & 0 & 0 & -\sqrt{t} & 0 & 0 & 0 & 0 & 0 & 0 & 0 \\
0 & 0 & 0 & 0 & 0 & 0 & 0 & 0 & 0 & -\sqrt{t} & 0 & 0 & 0 & 0 & 0 & 0 \\
0 & 0 & 0 & 0 & 0 & 0 & 0 & 0 & 0 & 0 & -\sqrt{t} & 0 & 0 & 0 & 0 & 0 \\
0 & 0 & 0 & 0 & 0 & 0 & 0 & 0 & 0 & 0 & 0 & -\sqrt{t} & 0 & 0 & 0 & 0 \\
0 & 0 & 0 & 0 & -\sqrt{t} & 0 & 0 & 0 & 1-t & 0 & 0 & 0 & 0 & 0 & 0 & 0 \\
0 & 0 & 0 & 0 & 0 & -\sqrt{t} & 0 & 0 & 0 & 1-t & 0 & 0 & 0 & 0 & 0 & 0 \\
0 & 0 & 0 & 0 & 0 & 0 & -\sqrt{t} & 0 & 0 & 0 & 1-t & 0 & 0 & 0 & 0 & 0 \\
0 & 0 & 0 & 0 & 0 & 0 & 0 & -\sqrt{t} & 0 & 0 & 0 & 1-t & 0 & 0 & 0 & 0 \\
0 & 0 & 0 & 0 & 0 & 0 & 0 & 0 & 0 & 0 & 0 & 0 & 1 & 0 & 0 & 0 \\
0 & 0 & 0 & 0 & 0 & 0 & 0 & 0 & 0 & 0 & 0 & 0 & 0 & 1 & 0 & 0 \\
0 & 0 & 0 & 0 & 0 & 0 & 0 & 0 & 0 & 0 & 0 & 0 & 0 & 0 & 1 & 0 \\
0 & 0 & 0 & 0 & 0 & 0 & 0 & 0 & 0 & 0 & 0 & 0 & 0 & 0 & 0 & 1
\end{bmatrix}$$

$$\sigma_1^{-1} = \begin{bmatrix}
1 & 0 & 0 & 0 & 0 & 0 & 0 & 0 & 0 & 0 & 0 & 0 & 0 & 0 & 0 & 0 \\
0 & 1 & 0 & 0 & 0 & 0 & 0 & 0 & 0 & 0 & 0 & 0 & 0 & 0 & 0 & 0 \\
0 & 0 & 1 & 0 & 0 & 0 & 0 & 0 & 0 & 0 & 0 & 0 & 0 & 0 & 0 & 0 \\
0 & 0 & 0 & 1 & 0 & 0 & 0 & 0 & 0 & 0 & 0 & 0 & 0 & 0 & 0 & 0 \\
0 & 0 & 0 & 0 & 1-\frac{1}{t} & 0 & 0 & 0 & -\frac{1}{\sqrt{t}} & 0 & 0 & 0 & 0 & 0 & 0 & 0 \\
0 & 0 & 0 & 0 & 0 & 1-\frac{1}{t} & 0 & 0 & 0 & -\frac{1}{\sqrt{t}} & 0 & 0 & 0 & 0 & 0 & 0 \\
0 & 0 & 0 & 0 & 0 & 0 & 1-\frac{1}{t} & 0 & 0 & 0 & -\frac{1}{\sqrt{t}} & 0 & 0 & 0 & 0 & 0 \\
0 & 0 & 0 & 0 & 0 & 0 & 0 & 1-\frac{1}{t} & 0 & 0 & 0 & -\frac{1}{\sqrt{t}} & 0 & 0 & 0 & 0 \\
0 & 0 & 0 & 0 & -\frac{1}{\sqrt{t}} & 0 & 0 & 0 & 0 & 0 & 0 & 0 & 0 & 0 & 0 & 0 \\
0 & 0 & 0 & 0 & 0 & -\frac{1}{\sqrt{t}} & 0 & 0 & 0 & 0 & 0 & 0 & 0 & 0 & 0 & 0 \\
0 & 0 & 0 & 0 & 0 & 0 & -\frac{1}{\sqrt{t}} & 0 & 0 & 0 & 0 & 0 & 0 & 0 & 0 & 0 \\
0 & 0 & 0 & 0 & 0 & 0 & 0 & -\frac{1}{\sqrt{t}} & 0 & 0 & 0 & 0 & 0 & 0 & 0 & 0 \\
0 & 0 & 0 & 0 & 0 & 0 & 0 & 0 & 0 & 0 & 0 & 0 & 1 & 0 & 0 & 0 \\
0 & 0 & 0 & 0 & 0 & 0 & 0 & 0 & 0 & 0 & 0 & 0 & 0 & 1 & 0 & 0 \\
0 & 0 & 0 & 0 & 0 & 0 & 0 & 0 & 0 & 0 & 0 & 0 & 0 & 0 & 1 & 0 \\
0 & 0 & 0 & 0 & 0 & 0 & 0 & 0 & 0 & 0 & 0 & 0 & 0 & 0 & 0 & 1
\end{bmatrix}$$

$$
\sigma_2 = \begin{bmatrix}
1 & 0 & 0 & 0 & 0 & 0 & 0 & 0 & 0 & 0 & 0 & 0 & 0 & 0 & 0 & 0 \\
0 & 1 & 0 & 0 & 0 & 0 & 0 & 0 & 0 & 0 & 0 & 0 & 0 & 0 & 0 & 0 \\
0 & 0 & 0 & 0 & -\sqrt{t} & 0 & 0 & 0 & 0 & 0 & 0 & 0 & 0 & 0 & 0 & 0 \\
0 & 0 & 0 & 0 & 0 & -\sqrt{t} & 0 & 0 & 0 & 0 & 0 & 0 & 0 & 0 & 0 & 0 \\
0 & 0 & -\sqrt{t} & 0 & 1-t & 0 & 0 & 0 & 0 & 0 & 0 & 0 & 0 & 0 & 0 & 0 \\
0 & 0 & 0 & -\sqrt{t} & 0 & 1-t & 0 & 0 & 0 & 0 & 0 & 0 & 0 & 0 & 0 & 0 \\
0 & 0 & 0 & 0 & 0 & 0 & 1 & 0 & 0 & 0 & 0 & 0 & 0 & 0 & 0 & 0 \\
0 & 0 & 0 & 0 & 0 & 0 & 0 & 1 & 0 & 0 & 0 & 0 & 0 & 0 & 0 & 0 \\
0 & 0 & 0 & 0 & 0 & 0 & 0 & 0 & 1 & 0 & 0 & 0 & 0 & 0 & 0 & 0 \\
0 & 0 & 0 & 0 & 0 & 0 & 0 & 0 & 0 & 1 & 0 & 0 & 0 & 0 & 0 & 0 \\
0 & 0 & 0 & 0 & 0 & 0 & 0 & 0 & 0 & 0 & 0 & 0 & -\sqrt{t} & 0 & 0 & 0 \\
0 & 0 & 0 & 0 & 0 & 0 & 0 & 0 & 0 & 0 & 0 & 0 & 0 & -\sqrt{t} & 0 & 0 \\
0 & 0 & 0 & 0 & 0 & 0 & 0 & 0 & 0 & 0 & -\sqrt{t} & 0 & 1-t & 0 & 0 & 0 \\
0 & 0 & 0 & 0 & 0 & 0 & 0 & 0 & 0 & 0 & 0 & -\sqrt{t} & 0 & 1-t & 0 & 0 \\
0 & 0 & 0 & 0 & 0 & 0 & 0 & 0 & 0 & 0 & 0 & 0 & 0 & 0 & 1 & 0 \\
0 & 0 & 0 & 0 & 0 & 0 & 0 & 0 & 0 & 0 & 0 & 0 & 0 & 0 & 0 & 1
\end{bmatrix}
$$

$$
\sigma_2^{-1} = \begin{bmatrix}
1 & 0 & 0 & 0 & 0 & 0 & 0 & 0 & 0 & 0 & 0 & 0 & 0 & 0 & 0 & 0 \\
0 & 1 & 0 & 0 & 0 & 0 & 0 & 0 & 0 & 0 & 0 & 0 & 0 & 0 & 0 & 0 \\
0 & 0 & 1-\frac{1}{t} & 0 & -\frac{1}{\sqrt{t}} & 0 & 0 & 0 & 0 & 0 & 0 & 0 & 0 & 0 & 0 & 0 \\
0 & 0 & 0 & 1-\frac{1}{t} & 0 & -\frac{1}{\sqrt{t}} & 0 & 0 & 0 & 0 & 0 & 0 & 0 & 0 & 0 & 0 \\
0 & 0 & -\frac{1}{\sqrt{t}} & 0 & 0 & 0 & 0 & 0 & 0 & 0 & 0 & 0 & 0 & 0 & 0 & 0 \\
0 & 0 & 0 & -\frac{1}{\sqrt{t}} & 0 & 0 & 0 & 0 & 0 & 0 & 0 & 0 & 0 & 0 & 0 & 0 \\
0 & 0 & 0 & 0 & 0 & 0 & 1 & 0 & 0 & 0 & 0 & 0 & 0 & 0 & 0 & 0 \\
0 & 0 & 0 & 0 & 0 & 0 & 0 & 1 & 0 & 0 & 0 & 0 & 0 & 0 & 0 & 0 \\
0 & 0 & 0 & 0 & 0 & 0 & 0 & 0 & 1 & 0 & 0 & 0 & 0 & 0 & 0 & 0 \\
0 & 0 & 0 & 0 & 0 & 0 & 0 & 0 & 0 & 1 & 0 & 0 & 0 & 0 & 0 & 0 \\
0 & 0 & 0 & 0 & 0 & 0 & 0 & 0 & 0 & 0 & 1-\frac{1}{t} & 0 & -\frac{1}{\sqrt{t}} & 0 & 0 & 0 \\
0 & 0 & 0 & 0 & 0 & 0 & 0 & 0 & 0 & 0 & 0 & 1-\frac{1}{t} & 0 & -\frac{1}{\sqrt{t}} & 0 & 0 \\
0 & 0 & 0 & 0 & 0 & 0 & 0 & 0 & 0 & 0 & -\frac{1}{\sqrt{t}} & 0 & 0 & 0 & 0 & 0 \\
0 & 0 & 0 & 0 & 0 & 0 & 0 & 0 & 0 & 0 & 0 & -\frac{1}{\sqrt{t}} & 0 & 0 & 0 & 0 \\
0 & 0 & 0 & 0 & 0 & 0 & 0 & 0 & 0 & 0 & 0 & 0 & 0 & 0 & 1 & 0 \\
0 & 0 & 0 & 0 & 0 & 0 & 0 & 0 & 0 & 0 & 0 & 0 & 0 & 0 & 0 & 1
\end{bmatrix}
$$

$$
\sigma_3 = \begin{bmatrix}
1 & 0 & 0 & 0 & 0 & 0 & 0 & 0 & 0 & 0 & 0 & 0 & 0 & 0 & 0 & 0 \\
0 & 0 & -\sqrt{t} & 0 & 0 & 0 & 0 & 0 & 0 & 0 & 0 & 0 & 0 & 0 & 0 & 0 \\
0 & -\sqrt{t} & 1-t & 0 & 0 & 0 & 0 & 0 & 0 & 0 & 0 & 0 & 0 & 0 & 0 & 0 \\
0 & 0 & 0 & 1 & 0 & 0 & 0 & 0 & 0 & 0 & 0 & 0 & 0 & 0 & 0 & 0 \\
0 & 0 & 0 & 0 & 1 & 0 & 0 & 0 & 0 & 0 & 0 & 0 & 0 & 0 & 0 & 0 \\
0 & 0 & 0 & 0 & 0 & 0 & -\sqrt{t} & 0 & 0 & 0 & 0 & 0 & 0 & 0 & 0 & 0 \\
0 & 0 & 0 & 0 & 0 & -\sqrt{t} & 1-t & 0 & 0 & 0 & 0 & 0 & 0 & 0 & 0 & 0 \\
0 & 0 & 0 & 0 & 0 & 0 & 0 & 1 & 0 & 0 & 0 & 0 & 0 & 0 & 0 & 0 \\
0 & 0 & 0 & 0 & 0 & 0 & 0 & 0 & 1 & 0 & 0 & 0 & 0 & 0 & 0 & 0 \\
0 & 0 & 0 & 0 & 0 & 0 & 0 & 0 & 0 & 0 & -\sqrt{t} & 0 & 0 & 0 & 0 & 0 \\
0 & 0 & 0 & 0 & 0 & 0 & 0 & 0 & 0 & -\sqrt{t} & 1-t & 0 & 0 & 0 & 0 & 0 \\
0 & 0 & 0 & 0 & 0 & 0 & 0 & 0 & 0 & 0 & 0 & 1 & 0 & 0 & 0 & 0 \\
0 & 0 & 0 & 0 & 0 & 0 & 0 & 0 & 0 & 0 & 0 & 0 & 1 & 0 & 0 & 0 \\
0 & 0 & 0 & 0 & 0 & 0 & 0 & 0 & 0 & 0 & 0 & 0 & 0 & 0 & -\sqrt{t} & 0 \\
0 & 0 & 0 & 0 & 0 & 0 & 0 & 0 & 0 & 0 & 0 & 0 & 0 & -\sqrt{t} & 1-t & 0 \\
0 & 0 & 0 & 0 & 0 & 0 & 0 & 0 & 0 & 0 & 0 & 0 & 0 & 0 & 0 & 1
\end{bmatrix}
$$

$$
\sigma_3^{-1} = \begin{bmatrix}
1 & 0 & 0 & 0 & 0 & 0 & 0 & 0 & 0 & 0 & 0 & 0 & 0 & 0 & 0 & 0 \\
0 & 1-\frac{1}{t} & -\frac{1}{\sqrt{t}} & 0 & 0 & 0 & 0 & 0 & 0 & 0 & 0 & 0 & 0 & 0 & 0 & 0 \\
0 & -\frac{1}{\sqrt{t}} & 0 & 0 & 0 & 0 & 0 & 0 & 0 & 0 & 0 & 0 & 0 & 0 & 0 & 0 \\
0 & 0 & 0 & 1 & 0 & 0 & 0 & 0 & 0 & 0 & 0 & 0 & 0 & 0 & 0 & 0 \\
0 & 0 & 0 & 0 & 1 & 0 & 0 & 0 & 0 & 0 & 0 & 0 & 0 & 0 & 0 & 0 \\
0 & 0 & 0 & 0 & 0 & 1-\frac{1}{t} & -\frac{1}{\sqrt{t}} & 0 & 0 & 0 & 0 & 0 & 0 & 0 & 0 & 0 \\
0 & 0 & 0 & 0 & 0 & -\frac{1}{\sqrt{t}} & 0 & 0 & 0 & 0 & 0 & 0 & 0 & 0 & 0 & 0 \\
0 & 0 & 0 & 0 & 0 & 0 & 0 & 1 & 0 & 0 & 0 & 0 & 0 & 0 & 0 & 0 \\
0 & 0 & 0 & 0 & 0 & 0 & 0 & 0 & 1 & 0 & 0 & 0 & 0 & 0 & 0 & 0 \\
0 & 0 & 0 & 0 & 0 & 0 & 0 & 0 & 0 & 1-\frac{1}{t} & -\frac{1}{\sqrt{t}} & 0 & 0 & 0 & 0 & 0 \\
0 & 0 & 0 & 0 & 0 & 0 & 0 & 0 & 0 & -\frac{1}{\sqrt{t}} & 0 & 0 & 0 & 0 & 0 & 0 \\
0 & 0 & 0 & 0 & 0 & 0 & 0 & 0 & 0 & 0 & 0 & 1 & 0 & 0 & 0 & 0 \\
0 & 0 & 0 & 0 & 0 & 0 & 0 & 0 & 0 & 0 & 0 & 0 & 1 & 0 & 0 & 0 \\
0 & 0 & 0 & 0 & 0 & 0 & 0 & 0 & 0 & 0 & 0 & 0 & 0 & 1-\frac{1}{t} & -\frac{1}{\sqrt{t}} & 0 \\
0 & 0 & 0 & 0 & 0 & 0 & 0 & 0 & 0 & 0 & 0 & 0 & 0 & -\frac{1}{\sqrt{t}} & 0 & 0 \\
0 & 0 & 0 & 0 & 0 & 0 & 0 & 0 & 0 & 0 & 0 & 0 & 0 & 0 & 0 & 1
\end{bmatrix}
$$

$$\mu^{\otimes 4} = \begin{bmatrix} 1 & 0 & 0 & 0 & 0 & 0 & 0 & 0 & 0 & 0 & 0 & 0 & 0 & 0 & 0 & 0 \\ 0 & t & 0 & 0 & 0 & 0 & 0 & 0 & 0 & 0 & 0 & 0 & 0 & 0 & 0 & 0 \\ 0 & 0 & t & 0 & 0 & 0 & 0 & 0 & 0 & 0 & 0 & 0 & 0 & 0 & 0 & 0 \\ 0 & 0 & 0 & t^2 & 0 & 0 & 0 & 0 & 0 & 0 & 0 & 0 & 0 & 0 & 0 & 0 \\ 0 & 0 & 0 & 0 & t & 0 & 0 & 0 & 0 & 0 & 0 & 0 & 0 & 0 & 0 & 0 \\ 0 & 0 & 0 & 0 & 0 & t^2 & 0 & 0 & 0 & 0 & 0 & 0 & 0 & 0 & 0 & 0 \\ 0 & 0 & 0 & 0 & 0 & 0 & t^2 & 0 & 0 & 0 & 0 & 0 & 0 & 0 & 0 & 0 \\ 0 & 0 & 0 & 0 & 0 & 0 & 0 & t^3 & 0 & 0 & 0 & 0 & 0 & 0 & 0 & 0 \\ 0 & 0 & 0 & 0 & 0 & 0 & 0 & 0 & t & 0 & 0 & 0 & 0 & 0 & 0 & 0 \\ 0 & 0 & 0 & 0 & 0 & 0 & 0 & 0 & 0 & t^2 & 0 & 0 & 0 & 0 & 0 & 0 \\ 0 & 0 & 0 & 0 & 0 & 0 & 0 & 0 & 0 & 0 & t^2 & 0 & 0 & 0 & 0 & 0 \\ 0 & 0 & 0 & 0 & 0 & 0 & 0 & 0 & 0 & 0 & 0 & t^3 & 0 & 0 & 0 & 0 \\ 0 & 0 & 0 & 0 & 0 & 0 & 0 & 0 & 0 & 0 & 0 & 0 & t^2 & 0 & 0 & 0 \\ 0 & 0 & 0 & 0 & 0 & 0 & 0 & 0 & 0 & 0 & 0 & 0 & 0 & t^3 & 0 & 0 \\ 0 & 0 & 0 & 0 & 0 & 0 & 0 & 0 & 0 & 0 & 0 & 0 & 0 & 0 & t^3 & 0 \\ 0 & 0 & 0 & 0 & 0 & 0 & 0 & 0 & 0 & 0 & 0 & 0 & 0 & 0 & 0 & t^4 \end{bmatrix}$$

The result comes out to be $-t^{-2}+t+t^2+t^4$, and since $\varepsilon(\beta)=-5$, we have $t^{\frac{1}{2}(\varepsilon(\beta)-n+1)}tr\left(\Phi_n(\beta)\mu^{\otimes n}\right)=\left(t^{-4}\right)\left(-t^{-2}+t+t^2+t^4\right)=-t^{-6}+t^{-3}+t^{-2}+1$. Hence the Jones polynomial for this knot is $-t^{-6}+t^{-5}-t^{-4}+2t^{-3}-t^{-2}+t^{-1}$.

Alice and Bob, having both computed this polynomial, can now form a secret key by agreeing upon some number derived from it. For instance, Alice and Bob might agree to take the value $\sum_i a_i\alpha_i^3$, where $a_i$ is the coefficient of the $i^{th}$ term and $\alpha_i$ is the exponent of $t$ in the $i^{th}$ term. In this case we would get a value of 108. (Clearly for more complicated knots this number might be quite a bit larger.) This value can now be used to confirm the identity of each party (i.e. they must agree otherwise we assume the other party to be Eve), as well as used to create a secret key.

# 3 Succession Models in Combinatorial Probability

Urn models have become a cornerstone of understanding combinatorial probability. Laplace used an urn model to show how the theories of Thomas Bayes can be used to infer the probability of successive events [40]. His was a model that given $k$ successes in $n$ trials (with replacement), we can have some idea of the probability of there being a success in the next trial. We present here a study of different models using Laplace's underlying theory and how they change given different processes at work. As well, we give a formula for the probability of succession given $n$ success in $n$ trials under an assumed binomial prior with and without replacements. (We also extend Laplace's formula for the case of an assumed uniform prior with replacements to that of an assumed prior without replacements.)

## 3.1. Introduction

The uses and abuses of probability have caused people to take strong positions on what should be considered probability and what should be deemed speculative and hence not a part of the mathematical theory. Using future events to gain insight into the past is one such area where mathematicians are divided (the two groups being known as the Frequentists and the Bayesians). What is meant here is the use of successive events to readjust one's assumptions as to what the underlying probability distribution might be. Given many such trials, the assumption of the probability distribution is assumed to get progressively more accurate. Underlying this theory, its foundation, is shaky ground. What is the underlying prior distribution? That is, what is our original assumption? What should it be? More importantly, Why? Given no real reasons why one distribution should be assumed as any better than any other, all results as to what the probability of some future event is becomes, at best, speculative. For this reason it is claimed that Laplace's model is inherently flawed. However, given some reason to assume one prior distribution over another and the model might become usable.

During the Second World War the Kengruppenbuch, the German cipher book, was captured by the British allowing cryptanalysts access to all the possible secret keys for the Enigma cipher – the code used to encrypt communications by the German military. The prevailing technique at that time used the relative frequency of each page of this book to create a prior distribution. Laplace's model was now at work in the world of code breaking.

Charged with breaking the Enigma cipher, Alan Turing and I. J. Good, took the prevailing theories and developed a new formula for estimating probabilities that did not lend itself to an intuitive understanding as to how (or why) it worked. Yet interestingly enough the Good-Turing estimator outperformed the more intuitive approaches. (See [49] for a recent look at this problem.)

The study of inference, and succession in particular, in classical probability therefore presents itself as a way to understand the underlying combinatorics of cryptological systems.

## 3.2. Succession in Trials with Replacement

### 3.2.1. Laplace and Uniform Distributions

An experiment is done. A ball is removed from an urn containing two balls of unknown color. You note the color of the ball just removed – it is white. You then return this ball to the urn and draw another ball. Again the ball is white and again you return the ball to the urn. What is the probability that if we would draw a ball from the urn once again, that it would be white?

This is the problem posed by Laplace [40]. His solution was as follows. There are two possible situations: either the urn has one white ball and one black ball (that is, a ball of some other color which we will call black) or it contains two white balls. The probability of drawing a white ball from the first urn is then ½, while the probability of drawing a white ball from the second urn is simply 1. The probability of drawing two white balls from the first urn (with replacement) is ¼, from the second urn it remains 1. Since each urn is equally likely, we can use Bayes's Theorem to find the probability that it is the first urn from which the balls were drawn and the probability that it was from the second urn. Bayes's Theorem states: $P(A \mid B) = \dfrac{P(B \mid A) \cdot P(A)}{P(B \mid A) \cdot P(A) + P(B \mid A^C) \cdot P(A^C)}$. Letting $A$ be the event that it is the first urn, $A^C$ the event that it is the second urn, and $B$ the event that two white balls were drawn, we have $P(A \mid B) = \dfrac{\frac{1}{4} \cdot \frac{1}{2}}{\frac{1}{4} \cdot \frac{1}{2} + 1 \cdot \frac{1}{2}} = \frac{1}{5}$. Similarly, the probability that the balls were drawn from the second urn comes out to be $\frac{4}{5}$. Using these values Laplace determined the probability that ‘the next ball drawn out of the urn would be white’ is $\frac{1}{2}\cdot\frac{1}{5}+1\cdot\frac{4}{5}=\frac{9}{10}$.

Laplace uses an assumption throughout his solution: Ignorance implies equal probability. Because we do not know the original distribution we can assume that the distribution was uniform. This idea can be well understood if we were to assume what we might call a ‘logical determiner.’ We use this phrase to indicate a general, combinatorial outcome: 3 white balls in an urn is simply 3 white balls in an urn. Thus if a person was to indicate the contents of an urn (and create the makeup), they would see 3 white balls as simply the fact that the urn contains three balls that are all white. The fact that they might have arisen under different distributions does not play a role. The only way we can assume this fact, that it does not play a role, is to assume that this ‘logical determiner’ played a role in the construction of the contents of the urn. To use Laplace’s famous rule of succession we then assume this logical determiner to be active.

Of interest to us is to see what happens as we pass to infinity. That is, when the number of balls in the urn increase without bound.

Let $X$ be the random variable indicating the total number of successes in $n$ trials, given repetitions. We can now calculate the probability of getting $k$ successes in $n$ trials as

$$P\{X=k \mid p\}=\binom{n}{k}p^k(1-p)^{n-k}\text{, for } k = 0, 1, 2,\ldots,n.$$

Letting $n=k$, and given the fact that the underlying probability, $p$, is uniformly distributed, we get:

$$P\{X=k\}=\frac{1}{G}\sum_{i=0}^{G}\left(\frac{i}{G}\right)^{k}=\frac{1}{G^{k+1}}\sum_{i=0}^{G}i^{k}\text{ , where } p=\left(\frac{i}{G}\right).$$

Using the formula for conditional probabilities, namely $P(A\mid B)=\frac{P(A\cap B)}{P(B)}$, the probability for success on the next try is:

$$P\{X=k+1\mid X=k\}=\frac{\frac{1}{G^{k+2}}\sum_{i=0}^{G}i^{k+1}}{\frac{1}{G^{k+1}}\sum_{i=0}^{G}i^{k}}=\frac{1}{G}\frac{\sum_{i=0}^{G}i^{k+1}}{\sum_{i=0}^{G}i^{k}}.$$

(This can also be expressed as $\frac{k+1}{k+2}\frac{B_{k+2}(G+1)-B_{k+2}}{B_{k+1}(G+1)-B_{k+1}}$. Where $B_p$ is the $p^{th}$ Bernoulli number, and $B_p(G)$ is the $p^{th}$ Bernoulli polynomial evaluated at *G*.)

Allowing $G\to\infty$ in our original equation (which in our model would mean that the total number of balls in the urn would get indefinitely large), we can pass to the integral

$$\int_0^1 P\{X=k\mid p\}dp=\int_0^1\binom{n}{k}p^k(1-p)^{n-k}dp=\binom{n}{k}\int_0^1 p^k(1-p)^{n-k}dp=\binom{n}{k}\frac{k!(n-k)!}{(n+1)!}=\frac{1}{n+1}.$$

(The identity $\int_0^1 p^k(1-p)^{n-k}dp=\frac{k!(n-k)!}{(n+1)!}$ can be shown by the use of the Beta function, $B(m+1,n+1)\equiv\frac{\Gamma(m+1)\Gamma(n+1)}{\Gamma(m+n+2)}=\frac{m!n!}{(m+n+1)!}=\int_0^1 u^m(1-u)^n du$, where $\Gamma(m)=(m-1)!$ is the Gamma function.)

Using Bayes's formula together with the formula for conditional probabilities, we can predict the next successive event.

$$P\{X = k+1 \mid X = k\} = \frac{\frac{k+1}{(n+1)(n+2)}}{\frac{1}{n+1}} = \frac{k+1}{n+2}$$

Setting $n = k$ we have $\frac{k+1}{k+2}$. This is Laplace's famous Rule of Succession. (The case in which $G$ remains constant but $k$ increases without bound is discussed in section 3.3.2.)

### 3.2.2. Binomial Distributions

What if this logical determiner is not available? That is, what if there is no reason for us to assume that the distribution is uniform. Our next logical step would be that instead of assuming that each possible combination has the same probability assignment we assume that each individual arrangement, including orderings, has the same probability. Let us take it as an assumption therefore that the distribution under consideration now is binomial. What would our formula be to predict the next outcome given that all the previous outcomes were successes?

Let us first give an example of the idea of predicting the next outcome given that all the previous trials (with replacements) resulted in successes where the underlying distribution of the probabilities is binomial. Instead of our logical determiner let us assume that we have the following scenario: An urn is to be filled with 3 balls. For each ball a fair coin is flipped: heads – a white ball goes in, tails – a black ball goes in. Not knowing the result of these flips, and thus the distribution of the white and black balls in the urn, we stick our hand into the urn and pull out a ball. It is white. We repeat this three times and find that each time a white ball is drawn. According to Laplace, which is when a uniform distribution of the probabilities is assumed, the probability of the next trial resulting in another white ball would be

$$P\{X=4 \mid X=3\}=\frac{P\{(X=4)and(X=3)\}}{P\{X=3\}}=\frac{\left(\frac{1}{3}\right)\cdot\left(\frac{1}{81}\right)+\left(\frac{1}{3}\right)\cdot\left(\frac{16}{81}\right)+\left(\frac{1}{3}\right)\cdot 1}{\left(\frac{1}{3}\right)\cdot\left(\frac{1}{27}\right)+\left(\frac{1}{3}\right)\cdot\left(\frac{8}{27}\right)+\left(\frac{1}{3}\right)\cdot 1}=\frac{\frac{98}{243}}{\frac{4}{9}}=\frac{49}{54}\approx .91$$

(The rule given earlier, that the probability is $\frac{n+1}{n+2}$, is only the limit as the number of balls gets larger.)

However, if we were to assume that the distributions of the probabilities was binomial, we would have

$$P\{X=4 \mid X=3\}=\frac{P\{(X=4)and(X=3)\}}{P\{X=3\}}=\frac{\left(\frac{3}{8}\right)\cdot\left(\frac{1}{81}\right)+\left(\frac{3}{8}\right)\cdot\left(\frac{16}{81}\right)+\left(\frac{1}{8}\right)\cdot 1}{\left(\frac{3}{8}\right)\cdot\left(\frac{1}{27}\right)+\left(\frac{3}{8}\right)\cdot\left(\frac{8}{27}\right)+\left(\frac{1}{8}\right)\cdot 1}=\frac{\frac{44}{216}}{\frac{1}{4}}=\frac{22}{27}\approx .81$$

The general form for the probability of success given that all the trials up until this point resulted in successes where the prior probabilities are binomially distributed over equal probabilities of success and failure (that is, we are assuming that the distribution of the initial probabilities is binomial with $p$ = ½ and $q$ = 1-$p$ = ½, i.e. each event in the Bernoulli trials is equally likely), is as follows:

Let $G$ be the number of balls in the urn (or the total number of elements in the sample space), $n$ the number of trials, and $k$ the number of successes. If $f(p)$ is the probability that $p$ is the prior probability (that is, the distribution in the urn), then we have

$$\sum P\{X=k \mid p\}f(p)=\sum_{i=0}^{G}\left[\binom{n}{k}p_i^{\,k}(1-p_i)^{n-k}\binom{G}{i}\frac{1}{2^G}\right], \text{ where } f(p)=\binom{G}{i}\frac{1}{2^G}.$$

Since $p_i = \frac{i}{G}$, we have $P\{X = k\} = \sum_{i=0}^{G}\left[\binom{n}{k}\left(\frac{i}{G}\right)^k\left(1-\frac{i}{G}\right)^{n-k}\binom{G}{i}\frac{1}{2^G}\right]$.

Letting $n = k$ we have the probability of $k$ successes in $k$ drawings (i.e. where every trial resulted in a success).

$$P\{X = k \mid n = k\} = \sum_{i=0}^{G}\left[\left(\frac{i}{G}\right)^k\binom{G}{i}\frac{1}{2^G}\right] = \frac{1}{2^G}\sum_{i=0}^{G}\left[\binom{G}{i}\left(\frac{i}{G}\right)^k\right] = \frac{1}{2^G G^k}\sum_{i=0}^{G}\left[\binom{G}{i}i^k\right].$$

The probability of succession using a Bayesian model over a binomial prior distribution, given $k$ successes in $k$ tries, is now:

$$P\{X = k+1 \mid X = k\} = \frac{\frac{1}{2^G G^{k+1}}\sum_{i=0}^{G}\binom{G}{i}i^{k+1}}{\frac{1}{2^G G^k}\sum_{i=0}^{G}\binom{G}{i}i^k} = \frac{1}{G}\frac{\sum_{i=0}^{G}\binom{G}{i}i^{k+1}}{\sum_{i=0}^{G}\binom{G}{i}i^k}. \qquad (*)$$

We can now use this to rework our earlier example. Given $G = 3$ and $n = k = 3$, we have

$$P\{X = 4 \mid X = 3\} = \frac{1}{3}\frac{\sum_{i=0}^{3}\binom{3}{i}i^4}{\sum_{i=0}^{3}\binom{3}{i}i^3} = \frac{22}{27} \approx .81$$

To take another example, let us assume the same situation as before, but this time with 5 balls in the urn instead of 3. According to Laplace we would have a probability of

$$P\{X=4 \mid X=3\}=\frac{P\{(X=4)and(X=3)\}}{P\{X=3\}}=\tfrac{1}{5}\frac{1^4+2^4+3^4+4^4+5^4}{1^3+2^3+3^3+4^3+5^3}=\frac{979}{1125}\approx .87$$

While under our binomial model we would have a probability of

$$P\{X=4 \mid X=3\}=\frac{P\{(X=4)and(X=3)\}}{P\{X=3\}}=\frac{1}{5}\frac{\binom{5}{1}1^4+\binom{5}{2}2^4+\binom{5}{3}3^4+\binom{5}{4}4^4+\binom{5}{5}5^4}{\binom{5}{1}1^3+\binom{5}{2}2^3+\binom{5}{3}3^3+\binom{5}{4}4^3+\binom{5}{5}5^3}=\frac{2880}{4000}=.72$$

As in the case of the uniform prior distribution we can look at what happens as the number of balls in the urn increases. Thus we are now trying to find the limiting case for the binomial distribution. To do this we use the following identity:

$$(1+x)^G=\binom{G}{0}+\binom{G}{1}x+\binom{G}{2}x^2+\cdots+\binom{G}{G}x^G.$$

Replacing $x$ with $e^x$, we have

$$(1+e^x)^G=\binom{G}{0}+\binom{G}{1}e^x+\binom{G}{2}e^{2x}+\cdots+\binom{G}{G}e^{Gx}.$$

Taking the $k^{th}$ derivative of this equation, the right side becomes

$$\binom{G}{0}0^k+\binom{G}{1}1^k e^x+\binom{G}{2}2^k e^{2x}+\cdots+\binom{G}{G}G^k e^{Gx}.$$

Letting $x = 0$, we have the summands in our expression (*). The left side of this equation, after taking the $k^{th}$ derivative, becomes $G(1+e^x)^{G-k} P(G^k)$, where $P(G^k)$ is a monic polynomial in $G$ of degree $k$. Letting $x = 0$ gives us $G2^{G-k} P(G^k)$. (*) now becomes:

$$\frac{\frac{1}{2^G G^{k+1}} G2^{G-k-1} P(G^{k+1})}{\frac{1}{2^G G^k} G2^{G-k} P(G^k)} = \frac{1}{2G} \frac{P(G^{k+1})}{P(G^k)} = \frac{1}{2} \frac{P(G^k)}{P(G^k)}.$$

Since the coefficients of the highest terms in the polynomials are both 1, as $G \to \infty$ this expression approaches $\frac{1}{2}$ (an answer which might seem quite anti-intuitive since $k$ can take on any finite value).

### 3.2.3. Applications

We have already discussed the uses of Laplace's model briefly. The binomial model, as well, lends itself to applications, examples of which are numerous especially in the theory of games. Imagine a gun in which each chamber is loaded via a process that gives each chamber a 50/50 chance of having (or not having) a bullet (perhaps a fair coin is tossed). Russian roulette is now played. (In Russian roulette the chambers of the gun are spun so that the chamber chosen is random.) If nobody has died so far (i.e. equivalent to the situation where there is replacement), what are the chances that the next player will find a bullet in the chamber? Or say we have a town with a very large population of homes with single occupancy, in which it just as equally likely that a man might live in any given home as there is that a woman lives there. Now let us assume that a computer randomly calls 400 homes with single occupancy. Given that the computer might call the same number multiple times, can we say what the probability that the computer will call a home in which a woman lives given the fact that all the calls up until this point have been to women? Given a large enough population we can say that it would be approximately 50/50. Or we can take as an example a board in which each square has an equal chance of having some event related to it as that of not having the event related to it, as, say, having a landmine. Given that on arbitrarily chosen squares (perhaps repeated) no landmine was found, what are the chances that given that so far none have blown up that one will on the next try?

The choice of which distribution should be assumed for any given scenario must be done carefully (unless, of course, it is simply given). For if it wasn't for the fact that we stated otherwise it would be more appropriate to assume a uniform distribution for the prior probabilities, as Laplace does, than a binomial model in our Russian roulette example. Since it is usually the number of bullets that is chosen (with equal probability), and not whether each individual chamber has a bullet.

## 3.3. Succession in Trials without Replacement

### 3.3.1. Uniform Distributions

We now extend the idea of a rule of succession to trials in which there are no replacements. Let us first work through an example by hand to see what we get. Assume we have an urn with five balls of at most two colors (where the different possible urns are uniformly distributed), and we draw three balls in a row without replacing them into the urn after each trial. If for each of the first three trials we pull out a white ball, what is the probability that on the fourth trial as well we will draw a white ball? There are three different possible urns (only three since the balls are not replaced after each trial): one containing 3 white balls, one containing 4 white balls, and one containing 5 white balls (i.e. all balls white). The probability of drawing three white balls in three trials would then be $\frac{1}{3} \cdot \frac{1}{10} + \frac{1}{3} \cdot \frac{4}{10} + \frac{1}{3} \cdot \frac{10}{10} = \frac{1}{2}$. The probability of drawing four white balls and three white balls (i.e. what simply results as the probability of drawing four white balls) in as many trials is $\frac{1}{3} \cdot 0 + \frac{1}{3} \cdot \frac{1}{5} + \frac{1}{3} \cdot \frac{5}{5} = \frac{2}{5}$. Dividing the latter value by the former we get the probability of drawing a fourth white ball given that three white balls were already drawn in as many trials. This comes out to be $\frac{4}{5}$.

Let us now look at the general case. To do this we will take a model much like Laplace's urn model, that is, with a uniform probability distribution for the prior probabilities, but assume that we do not have any replacements. We then have

$$P\{X=k\}=\frac{1}{G-k+1}\sum_{i=0}^{G}\frac{\binom{i}{k}\binom{G-i}{n-k}}{\binom{G}{n}},$$ where $\frac{1}{G-k+1}$ is the probability of any given distribution of balls in an urn. (We need not concern ourselves, or include, those urns that contain less than $k$ white balls, due to the fact that we are dealing with the case of no replacements. Therefore we use $\frac{1}{G-k+1}$.)

Allowing $n = k$ yields

$$\frac{1}{G-k+1}\frac{1}{\binom{G}{k}}\sum_{i=0}^{G}\binom{i}{k}=\frac{1}{G-k+1}\frac{1}{\binom{G}{k}}\binom{G+1}{k+1}=\frac{1}{G-k+1}\frac{G+1}{k+1}.$$

Given that $k$ white balls were drawn (without replacement), the probability that the $(k+1)^{\text{st}}$ ball will also be white, is $\frac{\frac{1}{G-k+1}\frac{G+1}{k+2}}{\frac{1}{G-k+1}\frac{G+1}{k+1}}=\frac{k+1}{k+2}$. ($\frac{1}{G-k+1}$ needs to be used for the numerator as well as for the denominator, since the probability is conditional, and the same conditions must apply to both.) (Notice how this formula yields the same result for our example as we got working it out by hand.) This result/formula is, interestingly enough, the exact same solution as when we had replacements. We sum up these points in a theorem.

**Theorem 3.1** Let $U$ be a collection (unordered) of $G$ objects of at most two different types (one of which we will deem the trial a success if selected), where the distribution of possible collections of objects is uniform. If an element is selected and discarded at random $k$ times ($k < G$) from $U$, with each trial resulting in a success, then the probability that the next selection will result in a success is $\frac{k+1}{k+2}$. If after each trial the selection was not discarded but rather returned, then the probability that the next trial will result in a success given that each of the $k$ trials resulted in success is $\frac{1}{G}\frac{\sum_{i=0}^{G} i^{k+1}}{\sum_{i=0}^{G} i^{k}}$.

Furthermore, as $G \to \infty$ this probability approaches the probability in the case in which there are no replacements.

What we have, therefore, is that our rule of succession for the model with replacements approaches that of the model without replacements as the size of the population increases. This makes a certain amount of intuitive sense, since the larger the collection the less the number of elements removed effects the distribution.

If, however, the size of our collection stays constant but more and more trials are performed, i.e. $k$ increases, then whether we were dealing with replacements or without replacements the probability would approach 1. Assuming, of course, that if no replacements are made $k$ remains less than $G$. This is seen as follows:

Rewriting $\frac{1}{G}\frac{\sum_{i=0}^{G} i^{k+1}}{\sum_{i=0}^{G} i^{k}}$, we get $\frac{\frac{1}{G}\left[1^{k+1} + 2^{k+1} + 3^{k+1} + \cdots + (G-1)^{k+1} + G^{k+1}\right]}{1^{k} + 2^{k} + 3^{k} + \cdots + (G-1)^{k} + G^{k}}$

$$= \frac{\left[1^k + 2^k + \cdots + G^k\right] - \left[1^k\left(\frac{G-1}{G}\right) + 2^k\left(\frac{G-2}{G}\right) + \cdots + (G-1)^k\left(\frac{1}{G}\right)\right]}{1^k + 2^k + 3^k + \cdots + (G-1)^k + G^k}$$

$$= 1 - \frac{1^k\left(\frac{G-1}{G}\right) + 2^k\left(\frac{G-2}{G}\right) + 3^k\left(\frac{G-3}{G}\right) + \cdots + (G-1)^k\left(\frac{1}{G}\right)}{1^k + 2^k + 3^k + \cdots + (G-1)^k + G^k}$$

$$> 1 - \frac{1^k + 2^k + 3^k + \cdots + (G-1)^k}{1^k + 2^k + 3^k + \cdots + (G-1)^k + G^k} > 1 - \frac{(G-1)^{k+1}}{G^k} = 1 - (G-1)\frac{(G-1)^k}{G^k}$$

$$= 1 - (G-1)\left(\frac{G-1}{G}\right)^k \to 1, \text{ as } k \to \infty .$$

Since $\frac{1}{G}\frac{\sum_{i=0}^{G} i^{k+1}}{\sum_{i=0}^{G} i^k}$ is bounded above by 1 (by definition of a probability space), we have

$$\frac{1}{G}\frac{\sum_{i=0}^{G} i^{k+1}}{\sum_{i=0}^{G} i^k} \to 1 \text{ as } k \to \infty .$$

## 3.3.2. Binomial Distributions

We now extend our ideas to the case in which the underlying distributions of the probabilities are not uniform, but rather binomial (with its underlying probability of success equal to its probability of failure, i.e. both ½, that is each event making up the Bernoulli trials is equally likely). We have in this case

$$P\{X=k\}=\sum_{i=0}^{G}\frac{\binom{i}{k}\binom{G-i}{n-k}}{\binom{G}{n}}\binom{G}{i}\frac{1}{2^G}=\frac{1}{2^G}\sum_{i=0}^{G}\frac{n!(G-n)!}{(n-k)!k!(i-k)!(G-i-n+k)!}.$$

(Where, as before, $\binom{G}{i}\frac{1}{2^G}$ is the probability of any given distribution of balls in an urn.)

Setting $n = k$ as we did before, we have

$$\frac{1}{2^G}\sum_{i=0}^{G}\frac{(G-k)!}{(i-k)!(G-i)!}=\frac{1}{2^G}\sum_{i=0}^{G}\binom{G-k}{i-k}=\frac{1}{2^G}\sum_{i=0}^{G-k}\binom{G-k}{i}=\frac{1}{2^G}2^{G-k}=\frac{1}{2^k}.$$

Note that $\sum_{i=0}^{G}\binom{G-k}{i-k}=\sum_{i=0}^{G-k}\binom{G-k}{i}$ since in our case we have only positive values for the binomial coefficients.

Our so-called Hypergeometric law of succession would therefore be: $\frac{\frac{1}{2^{k+1}}}{\frac{1}{2^k}}=\frac{1}{2}$. That is, the probability, of drawing a $(k+1)^{st}$ ball that is white from an urn (of which the distribution of the possible collections of balls is binomial, with probability of success for the Bernoulli trials being ½), given that $k$ balls have already be drawn (and discarded), all of which have been white, is ½. We sum up our results with the binomial distribution in the following theorem.

**Theorem 3.2** Let $U$ be a collection (unordered) of $G$ objects of at most two different types (one of which we will consider the trial a success if it is selected), where the distribution of possible collections of objects is binomial with its underlying probability of success equal to its probability of failure, i.e. both ½, that is each event making up the Bernoulli trials is equally likely. If an element is selected and discarded at random $k$ times ($k < G$) from $U$, with each trial resulting in a success, then the probability that the next selection will result in a success is ½. If after each trial the selection was not discarded but rather returned, then the probability that the next trial will result in a success given that each of the $k$ trials resulted in success is $\frac{1}{G}\frac{\sum_{i=0}^{G}\binom{G}{i}i^{k+1}}{\sum_{i=0}^{G}\binom{G}{i}i^{k}}$.

Furthermore, as $G \to \infty$ this probability approaches the probability of the case in which there are no replacements.

If the size of our collection, $G$, is kept constant, but more and more samples taken (with $k$ remaining less than $G$ when dealing without replacements), then under our model without replacements we would still have a probability of ½. However, in our model with replacements the probability approaches 1. This can be shown to be true as follows:

Rewriting $\frac{1}{G}\frac{\sum_{i=0}^{G}\binom{G}{i}i^{k+1}}{\sum_{i=0}^{G}\binom{G}{i}i^{k}}$, we get

$$\frac{\frac{1}{G}\left[\binom{G}{1}1^{k+1}+\binom{G}{2}2^{k+1}+\binom{G}{3}3^{k+1}+\cdots+\binom{G}{G-1}(G-1)^{k+1}+\binom{G}{G}G^{k+1}\right]}{\binom{G}{1}1^{k}+\binom{G}{2}2^{k}+\binom{G}{3}3^{k}+\cdots+\binom{G}{G-1}(G-1)^{k}+\binom{G}{G}G^{k}}$$

$$=\frac{\left[\binom{G}{1}1^{k}+\binom{G}{2}2^{k}+\cdots+\binom{G}{G}G^{k}\right]-\left[\binom{G}{1}1^{k}\left(\tfrac{G-1}{G}\right)+\binom{G}{2}2^{k}\left(\tfrac{G-2}{G}\right)+\cdots+\binom{G}{G-1}(G-1)^{k}\left(\tfrac{1}{G}\right)\right]}{\binom{G}{1}1^{k}+\binom{G}{2}2^{k}+\binom{G}{3}3^{k}+\cdots+\binom{G}{G-1}(G-1)^{k}+\binom{G}{G}G^{k}}$$

$$=1-\frac{\binom{G}{1}1^{k}\left(\tfrac{G-1}{G}\right)+\binom{G}{2}2^{k}\left(\tfrac{G-2}{G}\right)+\binom{G}{3}3^{k}\left(\tfrac{G-3}{G}\right)+\cdots+\binom{G}{G-1}(G-1)^{k}\left(\tfrac{1}{G}\right)}{\binom{G}{1}1^{k}+\binom{G}{2}2^{k}+\binom{G}{3}3^{k}+\cdots+\binom{G}{G-1}(G-1)^{k}+\binom{G}{G}G^{k}}$$

$>1-h(G)\dfrac{G(G-1)^{k}}{G^{k+1}}=1-h(G)\left(\dfrac{G-1}{G}\right)^{k}\to 1$, as $k\to\infty$, and where $h(G)$ is some (large) constant (in terms of *G*). And since $\dfrac{1}{G}\dfrac{\sum_{i=0}^{G}\binom{G}{i}i^{k+1}}{\sum_{i=0}^{G}\binom{G}{i}i^{k}}$ is bounded above by 1 (by definition of a probability space), we have $\dfrac{1}{G}\dfrac{\sum_{i=0}^{G}\binom{G}{i}i^{k+1}}{\sum_{i=0}^{G}\binom{G}{i}i^{k}}\to 1$ as $k\to\infty$.

(This is similar to how we showed earlier that $\dfrac{S_{k+1}(n)}{S_k(n)}\to n$ as $k\to\infty$, where $S_k(n)=\sum_{i=0}^{n}i^{k}$ .)

### 3.3.3. Applications

As examples go, we can easily see how some of our earlier examples can be viewed when no replacements are assumed. Given the computer that randomly calls homes with single occupancy in a certain town, it is easily seen that our newer model is more appropriate since it is natural to assume that once a number is called we would not have the computer call it again. Similarly with the case of the landmines, once a region containing a landmine is discovered that region would usually be excluded from all future investigations. Now if the area being searched is large enough (while the number of searches stays constant) we see that the probability of finding a landmine will be close to ½. So that given that we do not know the distribution, what came before does not change the chances in the next event.

## 3.4. Remarks

Underlying our remarks has been the notion of a second level distribution – a distribution of distributions. We have been attempting to understand the way in which the distribution of the possible distributions (of marbles in the urn) lies. Laplace naturally assumed that the (underlying) distribution was uniform; however, we showed that in many instances it is more reasonable to assume this distribution to be binomial. What has been underlying our search, therefore, has been the question 'What is the distribution of these distributions?' I.e. Is there any reason for us to assume it more likely that the underlying distribution is binomial rather than uniform? And if so, how might we quantify this notion so as to use it in our succession models? These questions and their relation to cryptography (prior assumptions and the possible distribution of those assumptions in data encryption) have not been the subject of much significant research.

# 4 Combinatorial Games on *n*-Simplex Graphs

## 4.1. Introduction

Blet is a game introduced by Sadun, Villegas, and Voloch [54]. The game consists of $2n$ coins placed in a circle with heads and tails alternating. The object of the game is to maximize the number of heads by turning over any head-tail-head triple to tail-head-tail, or turning any tail-head-tail triple to head-tail-head. In their paper the maximum value (number of heads) for any game is derived along with a number of additional problems that are raised. Its relation to the braid group $B_3$ is noticed as well. This relation between Blet and $B_3$ arises via the fact that the allowed moves in Blet are similar to the relator of the braid group $B_3$ under Artin's representation. (See section 2 for a full description of Braid Groups and Artin's representation of them.) The notion of having braid moves and relationships expressed as the flipping of coins motivated us to look at other games in which the relators were expressed as the flipping of coins.

The relation between graphs and games has been recognized throughout the development of Game Theory (see [10, 38, 61]). Graphs have played an integral role in the solution of coin flipping games, which are usually viewed as being upon graphs and grids. Coin flipping games, also known as $\sigma$-games, are games in which the coins adjacent to the one being flipped are also flipped (in certain versions only the adjacent coins are flipped). Coin flipping games are also popularly known as light-switching games. That is, a game in which a switch toggles neighboring switches as well. The objective in these games is, given an initial configuration, to turn all the lights off (or coins over). $\sigma$-games have been studied by Sutner [58, 59, 60] and by Barua and Ramakrishnan [9]. Anderson and Feil [2] studied the commercially available game "Lights Out," which is played on a $5 \times 5$ board. (For more on these games and their generalizations see [22, 27].) More recently, games have been looked at upon simplicial complexes [21]. Though any countable graph with vertices in n-dimensions can be viewed as a graph in two dimensions, the ability to view it in n-dimensions allows for the notion of motions and moves along the graph to have a more intuitive sense. We find simplicial complexes, therefore, to be a natural generalization. (See [62] for an earlier version of these ideas.)

We introduce a game in which each move, which we call a '*push*,' changes the labeling of every vertex in a (not necessarily maximal) clique (complete subgraph). Since non-maximal cliques can be viewed as n-simplexes, we will opt for the more visual notion of an n-simplex. There are some similarities to $\sigma$-games due to the label changing action of pushes; however, in many senses the similarity goes no further than the notion of a flipping game. In Sec. 4.5 we give a criterion to determine whether given two boards, which we define and call n-simplex graphs, it is possible to change one board into the other by a series of pushes. In Sec. 4.6 we determine (given that a solutions exists) how many different solutions there are. We use these results to find an algorithm, in Sec. 4.8, to determine a sufficient criterion for whether a graph is *n+1*-colorable, and in specific, whether a planar graph is 3-colorable.

While we will study here the game in all of its generality, the simple form of the game is one that is played upon a board of coins (or disks with each side a different color) that are tightly packed. That is, coins laid out so that the board is made up of little triangles of coins (i.e. each triangle consisting of three coins). A general board of this type can be formed by first taking three coins, each touching the other two, and then adding new coins, one at a time, so that each new coin touches at least two other coins (which each touch the other). The triangular and hexagonal boards are two examples of tightly packed boards (Fig. 4.1). A *push* would then consist of turning over any triangle that is made up of three coins (Fig. 4.2). The object of the game in its simplest form is to turn a board made up only of heads to a board made up only of tails. For small boards trial and error would be sufficient to determine whether and what solutions exist. And given enough time, for boards of specific shapes, certain patterns might arise with regards to the solutions.

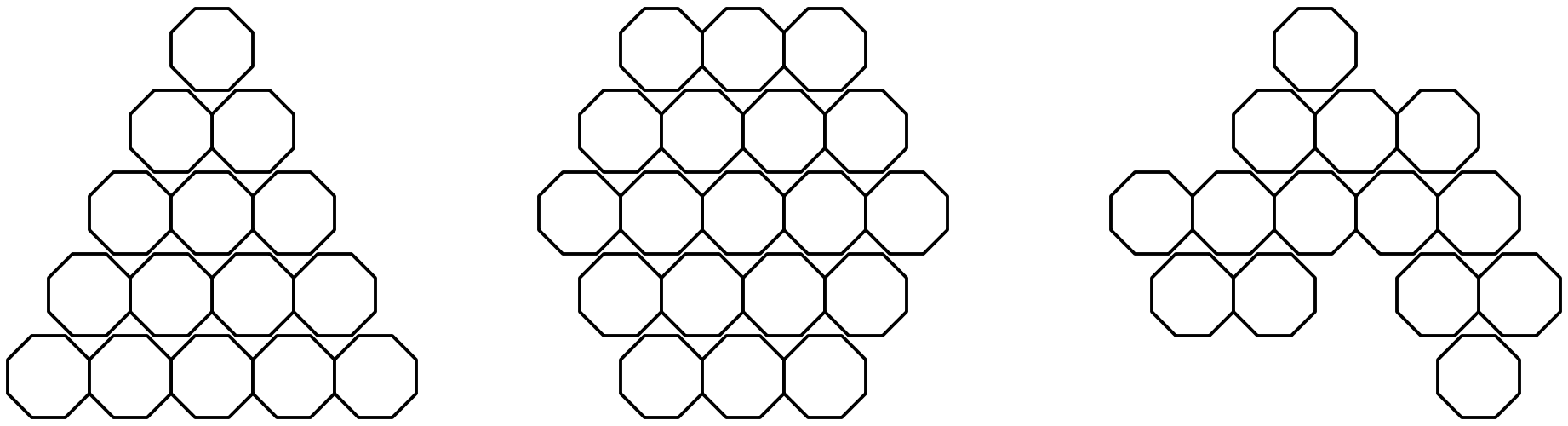

Figure 4.1. A triangular, hexagonal, and arbitrarily shaped board.

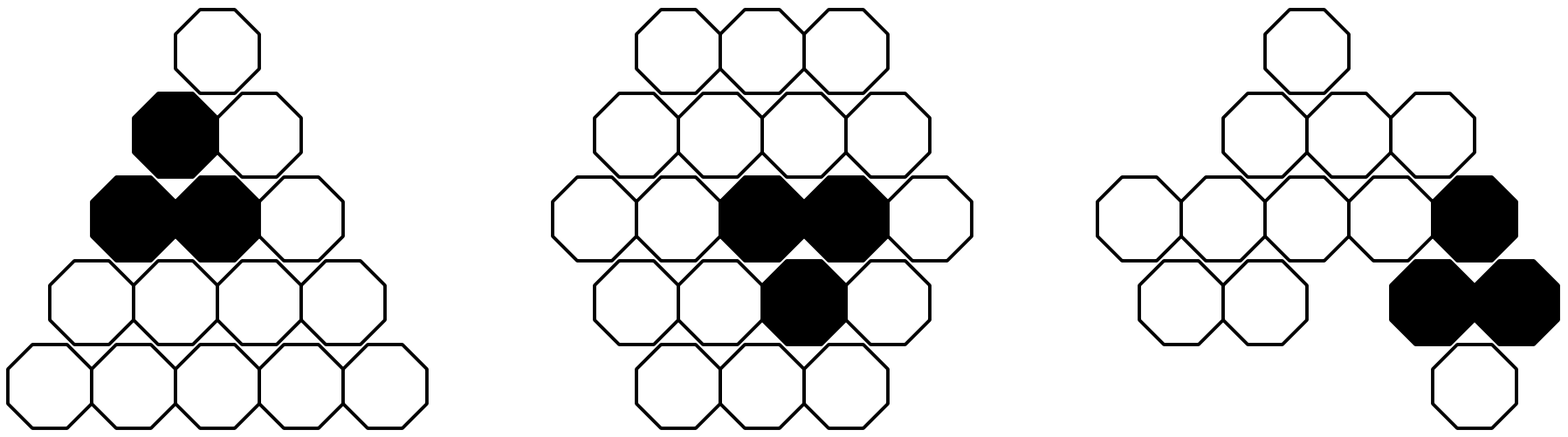

Figure 4.2. A few different pushes.

Though in its simplest form the game is played with the objective of turning a board made up exclusively of disks (coins) of one color (heads) to disks of the alternate color (tails), the game can be generalized to any two colorings of a board. That is, given two colorings of a board, where a coloring is simply the assignment of different colors to different disks, we call one board the initial board and the other board the terminal board. The object of the game then, is to find a collection of pushes (a 'solution set') that turn the coloring of the initial board into the coloring of the terminal board. We mention 'collection,' as opposed to an ordered set, because, as is easily seen, pushes are commutative, and hence the order of the pushes in the solution set is unimportant.

## 4.2. Definitions

Following Hatcher [28], we first define an n-simplex and a simplicial complex topologically and then indicate their natural corresponding graphical interpretation.

**Definition 4.2.1** *Given any set* $V = \{v_0, v_1, ..., v_n\}$ *of* $n+1$ *points in* $\mathbf{R}^N$*, such that the differences* $v_1 - v_0, v_2 - v_0, ..., v_n - v_0$ *are linearly independent, the* n-simplex with vertices V *is the convex hull of V, i.e. the set of all points of the form* $t_0 v_0 + t_1 v_1 + ... + t_n v_n$*, where* $\sum_{i=0}^{n} t_i = 1$ *and* $t_i \geq 0$ *for all i.*

**Definition 4.2.2** *A* simplicial complex $\Delta$ on a finite set *V is a collection of subsets of V such that*

*i)* $\{v\} \in \Delta$ *for all* $v \in V$.

*ii) if* $F \in \Delta$ *and* $G \subseteq F$ *then* $G \in \Delta$.

*The members of* $\Delta$ *are called* simplices *or* faces, *and the elements of V are called vertices.*

The use (and advantage) of topological notions is that it will allow us to define n-simplex graphs upon any compact connected *q*-manifold, with $q \geq n$.

Viewed purely from a graph theoretic perspective, an n-simplex is simply a complete graph, $U_i$, such that $|V(U_i)| = n+1$. Thus the set of vertices, $V(U_i)$, uniquely determines the n-simplex. And a simplicial complex is merely a graph $G = \bigcup U_j$ with each $U_j$ complete. Due to the natural mapping from the topological definition of simplicial complexes to the graph-theoretic definition, we will use them interchangeably, adopting the notion of an n-dimensional graph, *G*, comprised of n-simplexes. It should be noted that the connection to the topological notion has certain limitations since we are not concerned with the topology of n-simplexes or with viewing them as convex hulls, but rather with their vertex sets, edge sets, and faces. Therefore, except for the intuition that it brings with it, the ideas are primarily graph-theoretic, and should be viewed as such.

We now define the notion of an n-simplex graph based upon the notions of a simplicial complex.

**Definition 4.2.3** *G is said to be an* n-simplex graph *if:*

*i)* *G is a simplicial complex.*

*ii)* *For any* $F \in G$ *with* $|F| < n+1$ *there exists a* $K \in G$ *such that* $F \subset K$.

Alternately an n-simplex *graph* can be defined (and viewed) as a graph $G = \bigcup U_j$ where each $U_j$ is complete and $|V(U_j)| = n+1$.

We take it as a general assumption throughout this paper that all graphs are of finite size.

## 4.3. Region-Paths and Region-Connected Graphs

The building blocks of n-simplex graphs are n-simplexes. Topologically they are the convex hulls of the set of vertices. It is therefore natural that we shall refer to these n-simplexes as the regions of *G*.

We now extend the familiar notions of adjacent vertices and paths along vertices to that of regions. Adjacent regions are defined in the same vein as that of adjacent vertices.

**Definition 4.3.1** *Two n-simplexes,* $S_i$ *and* $S_j$ *are said to be* adjacent *if* $S_i \cap S_j$ *is an (n-1)-simplex.*

**Definition 4.3.2** *A* region-path *from* $S_i$ *to* $S_j$ *is a set of n-simplexes* $\{S_i, S_{i+1}, \ldots, S_j\}$ *such that* $S_k$ *is adjacent to* $S_{k+1}$ *for* $i \le k < j$.

Let $L(G): V(G) \to Z_n$ be a labeling of the vertices of a graph $G$ from the set $\{0,1,\ldots,n-1\}$. We call a move a *push* if it is a function $f_{S_i}: L_1(G) \to L_2(G)$, acting on an n-simplex $S_i = \{v_0, v_1, \ldots, v_n\}$, such that $f_{S_i}[L_1(v_j)] = L_1(v_j) + 1 \pmod{n}$ for $v_j \in S$.

**Definition 4.3.3** *A graph is said to be* region-connected *if given any two n-simplexes* $S_1, S_2 \in G$*, there exists a region-path connecting them.*

Pictorially, a *region-connected* n-simplex graph can be viewed as a graph formed by pasting n-simplexes together by their ($n$-1)-simplexes. That is, the intersection of two n-simplexes that are pasted together is an ($n$-1)-simplex. For example, we could form a region-connected 2-simplex graph by pasting triangles together by their edges, or form a region-connected 3-simplex graph by pasting tetrahedrons together by their triangles.

## 4.4. An Invariant Under Pushes

We pose the following two general problems concerning motions, or re-labelings, of graphs:

**Question 4.4.1** *Let G be any region-connected n-simplex graph. Given two labelings,* $L_1(G)$ & $L_2(G)$*, does there exist a series of pushes with which we can change* $L_1(G)$ *into* $L_2(G)$*?*

**Question 4.4.2** *Given that a solution to Question 4.4.1 exists, how many different solutions are there?*

We leave the solution to Question 4.4.2 until Sec. 4.6.

We wish to find a class of n-simplex graphs under which our question is solvable. We claim that if *G* is a region-connected n-simplex graph then the question is solvable (not that it is affirmative) if $\chi(G) = n+1$. We start by finding a function of $L(G)$ which is invariant under pushes.

Let $Z_m$ be the labeling set for $L_1$ & $L_2$. ($Z_m$ can easily be allowed to be the larger of the two labeling sets if they differ.) Since $\chi(G) = n+1$, *G* can be colored with the set $\{i_0, i_1, ..., i_n\}$, where $i_k$ is as follows:

$$i_k = \begin{bmatrix} I_k & & & 0 & 0 \\ & \ddots & & & 0 \\ & & e^{\frac{i2\pi}{m}} & & \\ 0 & & & \ddots & \\ 0 & 0 & & & I_{n-k-1} \end{bmatrix} \quad 0 \le k < n,\ i_n = \begin{bmatrix} e^{\frac{i2(m-1)\pi}{m}} & & 0 \\ & \ddots & \\ 0 & & e^{\frac{i2(m-1)\pi}{m}} \end{bmatrix} = e^{\frac{i2(m-1)\pi}{m}} I_n$$

We have $i_0^m = i_1^m = \cdots = i_{n-1}^m = i_n^m = i_0 i_1 \cdots i_n = I$.

Given that $\{i_0, i_1, ..., i_{n-1}\}$ is linearly independent and that the order of each $i_j$ is *m*, we find that $\langle i_0, i_1, ..., i_{n-1} \rangle \approx \underbrace{Z_m \times Z_m \times \ldots \times Z_m}_{n-times}$.

Let $i(v_j) \in \{i_0, i_1, \cdots, i_n\}$ be the coloring of the vertex $v_j \in V(G)$ and $l(v_j) \in L(G)$ its label. Assign the value $i^{l(v_j)}(v_j)$ to the vertex $v_j$, for each $v_j \in V(G)$.

Let $P[L(G)] = \prod_{v_j \in V(G)} i^{l(v_j)}(v_j)$.

**Lemma 4.4.3** *$P[L(G)]$ is invariant under pushes.*

**Proof** Let $\{i_0^{\alpha_0}, i_1^{\alpha_1}, \cdots, i_n^{\alpha_n}\}$ be the set of values assigned to the vertices of an arbitrary n-simplex in *G*. A push on this n-simplex would send $i_j^{\alpha_j} \mapsto i_j^{\alpha_j + 1}$ $\forall j \in \{0,1,...,n\}$. Thus we have:

$$i_0^{\alpha_0} i_1^{\alpha_1} \cdots i_n^{\alpha_n} \rightarrow i_0^{\alpha_0+1} i_1^{\alpha_1+1} \cdots i_n^{\alpha_n+1} = i_0^{\alpha_0} i_1^{\alpha_1} \cdots i_n^{\alpha_n} \left(i_0 i_1 \cdots i_n\right) = i_0^{\alpha_0} i_1^{\alpha_1} \cdots i_n^{\alpha_n} \quad \text{(since } i_0 i_1 \cdots i_n = I \text{)}.$$

(Since the push acts only upon a given n-simplex, the remainder of $G$ remains unchanged.)

$P[L(G)]$ is therefore invariant under pushes.

## 4.5. Solutions for a Specific Class of Graphs

We now find a class of graphs under which our question is solvable and present an algorithm for finding a series of pushes given that a solution exists.

**Theorem 4.5.1** *Let G be a region-connected n-simplex graph with* $\chi(G) = n + 1$. *Then there exists a set F, of pushes, such that* $F[L_1(G)] = L_2(G)$ *iff* $P[L_1(G)] = P[L_2(G)]$.

**Proof** By Lemma 4.4.3 $P[L(G)]$ is invariant under pushes; the necessity of the equality therefore follows. (The condition of region-connectedness is not actually required for this direction.)

To show the sufficiency of the equality we will construct a solution. Since $G$ is region-connected, we can find $n$+1 region-paths in G such that:

(i) $s_{1,i_j}, s_{2,i_j}, \cdots, s_{t,i_j}$ is a region-path of n-simplexes with $s_{(2k-1),i_j} = \{v_0, v_1, \cdots, v_n\}$, $s_{(2k),i_j} = \{v_1, \cdots, v_n, v_{n+1}\}$, then $s_{(2k-1),i_j} \cap s_{(2k),i_j} = \{v_1, \cdots, v_n\}$, $k = \{1, 2, \ldots\}$, where $i_j$ is the coloring of both $v_0$ and $v_{n+1}$.

(ii) Given any $v \in V(G)$ with coloring $i_j$, $v$ is in $s_{r,i_j}$ for some $r$.

(iii) $s_{t,i_0} = s_{t,i_1} = \cdots = s_{t,i_n}$. I.e. all paths end with the same n-simplex.

(No other conditions exist for these paths; thus a region-path might cross or retrace itself.)

Let us now assume that two labelings, $L_1$ & $L_2$, differ in only one n-simplex. Let $f_{1,i_0}^{k_{1,i_0}}$ be a push to such a power that the vertex in this n-simplex colored $i_0$ attains the same labeling as in $L_2(G)$. We claim that $L_1(G) = L_2(G)$. Since $P[L_1(G)] = P[L_2(G)]$ (by assumption) we have $i_0^{\alpha_{1,0}} i_1^{\alpha_{1,1}} \cdots i_n^{\alpha_{1,n}} = i_0^{\alpha_{2,0}} i_1^{\alpha_{2,1}} \cdots i_n^{\alpha_{2,n}}$. And given that $i_0^{\alpha_{1,0}} = i_0^{\alpha_{2,0}}$, we have $A_{ij} = i_1^{\alpha_{1,1}}, \cdots, i_n^{\alpha_{1,n}} = i_1^{\alpha_{2,1}}, \cdots, i_n^{\alpha_{2,n}} = B_{ij}$. As well, since $a_{1,1} = e^{\frac{i2\alpha_{1,n}(m-1)\pi}{m}} = e^{\frac{i2\alpha_{2,n}(m-1)\pi}{m}} = b_{1,1}$, we have $\alpha_{1,n} = \alpha_{2,n}$. That the remaining labels are also equal can be seen as a result of their linear independence. We therefore have $\alpha_{1,j} = \alpha_{2,j} \forall j \in \{0,1,\cdots,n\}$.

Let us now assume that the labelings differ arbitrarily. We form $n+1$ sequences of pushes as follows: For the coloring $i_j$ we have the sequence $f_{1,i_j}^{p_{1,i_j}} f_{2,i_j}^{m-p_{1,i_j}} f_{3,i_j}^{p_{3,i_j}} f_{4,i_j}^{m-p_{3,i_j}} \cdots$, where the pushes act upon the n-simplexes of the region-paths given above, and where the $(2k-1)^{st}$ term is of the form $f_{(2k-1),i_j}^{p_{(2k-1),i_j}}$ and the $(2k)^{th}$ term is of the form $f_{(2k),i_j}^{m-p_{1,i_j}}$, for $k=\{1,2,\ldots\}$, and where $p_{(2k-1),i_j}$ is the power necessary so that if $v$ is the vertex in the n-simplex colored $i_j$, then $l_1\left[f_{(2k-1),i_j}^{p_{(2k-1),i_j}}(v)\right]=l_2(v)$. I.e. the labeling of that vertex is the same as in $L_2(G)$.

The last term in each of the sequences is either of the form $f_{(2k-1),i_j}^{p_{(2k-1),i_j}}$ or $f_{(2k),i_j}^{m-p_{(2k-1),i_j}}$, depending upon whether there is an odd or an even number of elements in the region-path.

Note that since $s_{(2k-1),i_j} \cap s_{(2k),i_j} = \{v_1,\cdots,v_n\}$, $f_{(2k-1),i_j}^{p_{(2k-1),i_j}} f_{(2k),i_j}^{m-p_{(2k-1),i_j}}$ raises the values of all those vertices not colored $i_j$ by *m*. Therefore, all vertices, except those colored $i_j$, remain unchanged.

Each sequence will therefore change the labelings of all those vertices of *G* colored $i_j$, except, possibly, for that of the last n-simplex of the path (upon which the sequence acts). This being true for each $i_j$, we need only concern ourselves now with this last n-simplex, which, by construction, is the same for each path (the rest of *G* can therefore be ignored). Having already proven the theorem for that case our proof is complete.

## 4.6. Board Classes and Solution Sizes

Our objective in this section is to show how many different solutions exist when there is a solution. To accomplish this we show that the set of labelings of *G* form equivalence classes. This we do because we wish to demonstrate that the sizes of all these equivalence classes are the same.

Let us add to our collection of pushes acting on *G*, *n* additional elements acting only upon one vertex out of the $n+1$ vertices in some given n-simplex. Without loss of generality we can assume that these new elements all act upon the same n-simplex $S_0 = \{v_0, v_1, ..., v_n\}$. We now let our *n* new moves be such that $h_i[L_1(v_i)] = L_1(v_i) + 1$ (mod *n*) for only the given $v_i \in S$ with $i = \{0, 1, \ldots, n-1\}$. (In other words, these *n* new pushes will only change one vertex from the given n-simplex, all the other vertices will remain unchanged.) Since $P[L(G)]$ can now take on any value (due to the addition of these new moves), it is possible to change from any given labeling of *G* to any other labeling of *G*.

We now form classes of labelings of *G* such that two labelings, *K* and *L*, are in the same class if $P[K(G)] = P[L(G)]$.

**Definition 4.6.1** *Let $L_1(G)$ & $L_2(G)$ be two labelings of a region-connected n-simplex graph G, with the condition that $\chi(G) = n+1$. Then we say that $L_1(G)$ is* label-equivalent *to $L_2(G)$, written $L_1(G) \sim L_2(G)$, if $P[L_1(G)] = P[L_2(G)]$.*

**Lemma 4.6.2** *The relation* $L_1(G) \sim L_2(G)$ *is an equivalence relation.*

**Proof** Using Theorem 4.5.1, it can be shown that the conditions for an equivalence relation are easily satisfied.

Since an equivalence relation decomposes a set into mutually disjoint subsets (see [29] for a simple proof), we have the following:

**Corollary 4.6.3** *The equivalence relation ~ provides a decomposition of the set of all labelings of G into distinct (mutually disjoint) equivalence classes.*

Let *m* be the size of the labeling set of *G* and *n* be such that $\chi(G) = n + 1$. We now have the further fact, that

**Corollary 4.6.4** *There are* $m^n$ *distinct equivalence classes.*

**Proof** The proof follows easily from Theorem 4.5.1 and the fact that $P[L(G)]$ can be anyone of $m^n$ different values.

Let $h_i^{\alpha_i}$ be a mapping from one class of labelings to another class, then if *K* and *L* are two labelings in some class, $h_i^{\alpha_i} K$ and $h_i^{\alpha_i} L$ will be two labelings in this other class.

However, it is easily shown that $h_i^{\alpha_i} K = h_i^{\alpha_i} L$ if and only if $K = L$. Therefore there is a one to one relationship between the different equivalence classes, and we have

**Corollary 4.6.5** *The equivalence relation ~ divides the set of labeled graphs into equivalence classes of equal size.*

By Corollary 4.6.4 we know that there are $m^n$ different values that *P* can take. Letting $|V(G)| = v$, there are $m^v$ different labelings of *G*. Since, by Corollary 4.6.5, the sizes of the equivalence classes are the same, there are $m^{v-n}$ different labelings for each class. Therefore, we have that

**Corollary 4.6.6** *There are exactly $m^{v-n}$ elements in each equivalence class formed by the equivalence relation ~.*

We now form classes of pushes that all act in a similar manner on a labeling. Let $R(G)$ be the set of all n-simplexes (i.e. regions) in *G*.

**Definition 4.6.7** *Let f and g be two words in $R(G)$ and L some labeling of G, then we say that f is* congruent *to g, written as f ≡ g, if* $f[L(G)] = g[L(G)]$.

**Lemma 4.6.8** *The relation f ≡ g is an equivalence relation.*

**Proof** Follows from the definition of an equivalence relation.

As before, since an equivalence relation creates a decomposition of the set into mutually disjoint subsets (see [29]), we have the following:

**Corollary 4.6.9** *The equivalence relation $\equiv$ provides a decomposition of the set of all words in $R(G)$ into distinct (mutually disjoint) equivalence classes.*

**Corollary 4.6.10** *The equivalence relation $\equiv$ divides $R(G)$ into equivalence classes of equal size.*

**Proof** We form a mapping $\varphi : K \to L$ from one equivalence class, *K*, into another equivalence class, *L*, by $f \mapsto \mu f$, where $\mu$ is such that $\mu f \in L$. It can easily be shown that $\mu f \equiv \mu g$ iff $f \equiv g$, and $\mu f = \mu g$ iff $f = g$. Therefore, since $K$ and $L$ were arbitrary, there is a one to one correspondence between equivalence classes.

Let $|R(G)| = r$, then there are $m^r$ different sets of moves possible on *G*. Since we know by Corollary 4.6.10 that the classes of moves are each the same size, and by Corollary 4.6.6 that there are $m^{v-n}$ different classes of labelings upon which these moves act, we have that there are $\frac{m^r}{m^{v-n}} = m^{r-v+n}$ different sets of pushes for each class of labelings. Thus we have proved

**Theorem 4.6.10** *Given a graph G, as in Theorem 4.5.1, the number of solutions that exist with which one labeling can be changed into another labeling, so long as a solution exists, is* $m^{r-v+n}$.

## 4.7. Examples of Games on 2-Simplexes

**Example 4.7.1** As an example, imagine the board discussed in section 4.1. The board is made up of coins (or disks with each side a different color), which are tightly packed. That is, coins laid out so that the board is made up of little triangles of coins (i.e. each triangle consisting of three coins). A general board of this type can be formed by first taking three coins, each touching the other two, then adding new coins, one at a time, so that each new coin touches at least two other coins (which each touch the other). For example, the board can be shaped hexagonally or triangularly (each row having one more coin than the previous row). Our question would now be stated as follows: Given any two of these boards (and an arrangement of the coins in heads and tails for each), does there exist a set of pushes such that one board can be changed into the other?

Since $\chi(G) = 3$ (which can be easily demonstrated) a solution is now easily determined (with the existence of a solution depending on whether or not $P[L_1(G)] = P[L_2(G)]$).

For example, say we would like to change 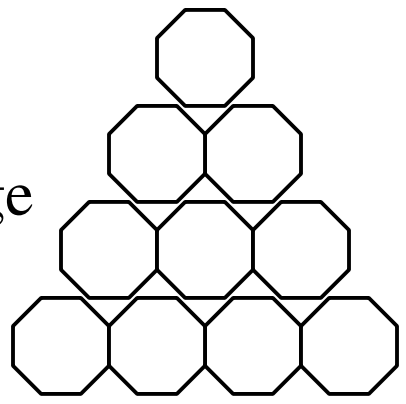 → 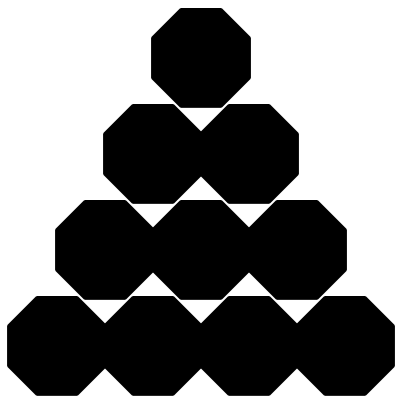

Then since the $P$ values of both of these boards are not equal (i.e. $P[L_1(G)] \neq P[L_2(G)]$, for any 3-coloring of it) no solution exists.

However, if we are given the following two boards: 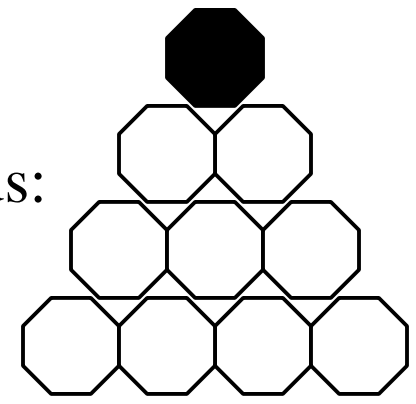 → 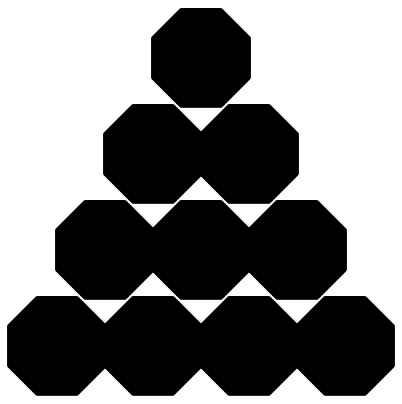

then since the values are equal (i.e. $P[L_1(G)] = P[L_2(G)]$, if we would 3-color them), there does exists a solution, and since $m^{r-v+n} = 2^{9-10+2} = 2$, we find that there are in fact two.

We can visualize these solutions as follows: If we replace each disk with a vertex, and for every two disks that touch assign an edge, then we can associate with each board an underlying graph, $G$ (Fig. 4.3). The underlying graph of a board can be thought of as little triangles pasted together by their edges. We now associate with each color of a disk a labeling of its respective vertex. Choosing one side of a disk as 'face down' and the other side as 'face up' we might label a vertex 0 if it was face down and 1 if it was face up. Therefore, based on the coloring of the board we attain a labeling of the graph. Due to this association we will often refer to the board as the graph and to the graph as the board.

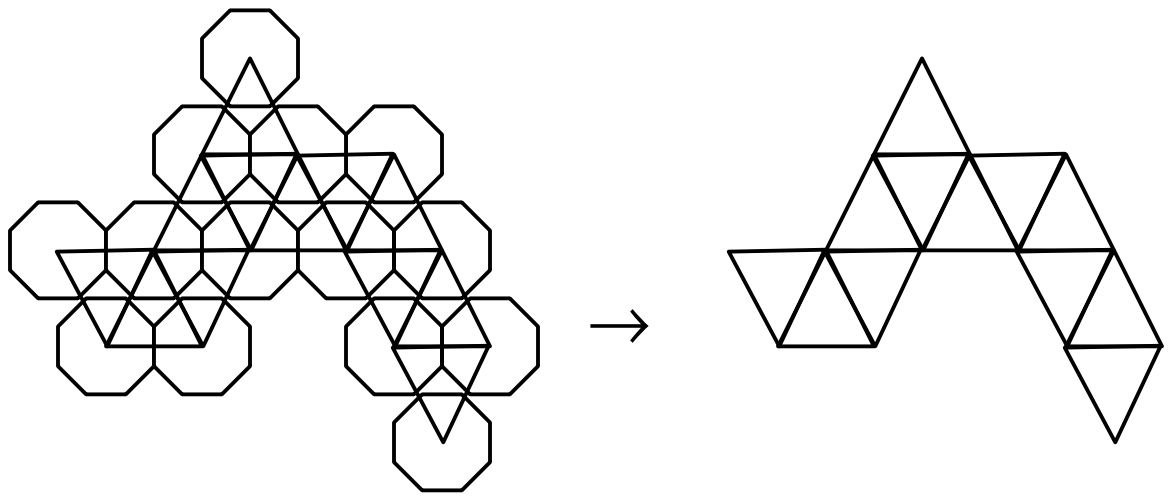

Figure 4.3. The underlying graph of a board

The solutions to our example can now been seen to be

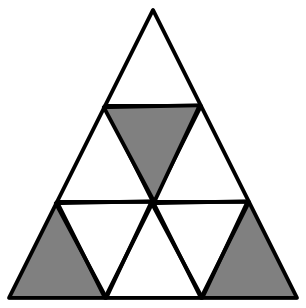 and 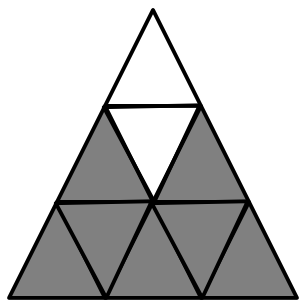

(Where 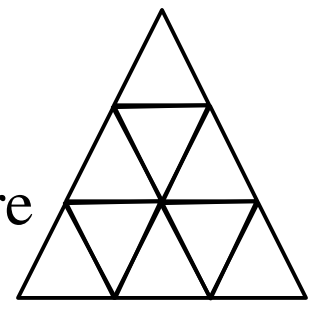 is obtained by connecting the disks as such 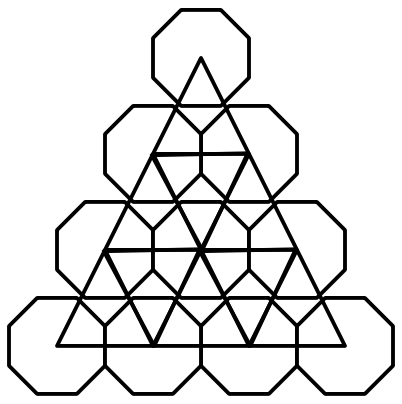

The vertices of the colored triangles are then flipped, that is, the pushes act upon the colored triangles).

If we would now increase our board size to five rows, and have as our initial and terminal boards all face down and all face up disks (respectively), then not only would the *P* values be equal, and hence a solution would exist, but we would find that since $2^{r-v+2} = 2^{16-15+2} = 2^3$, there would be eight different solutions. (See Fig. 4.4.)

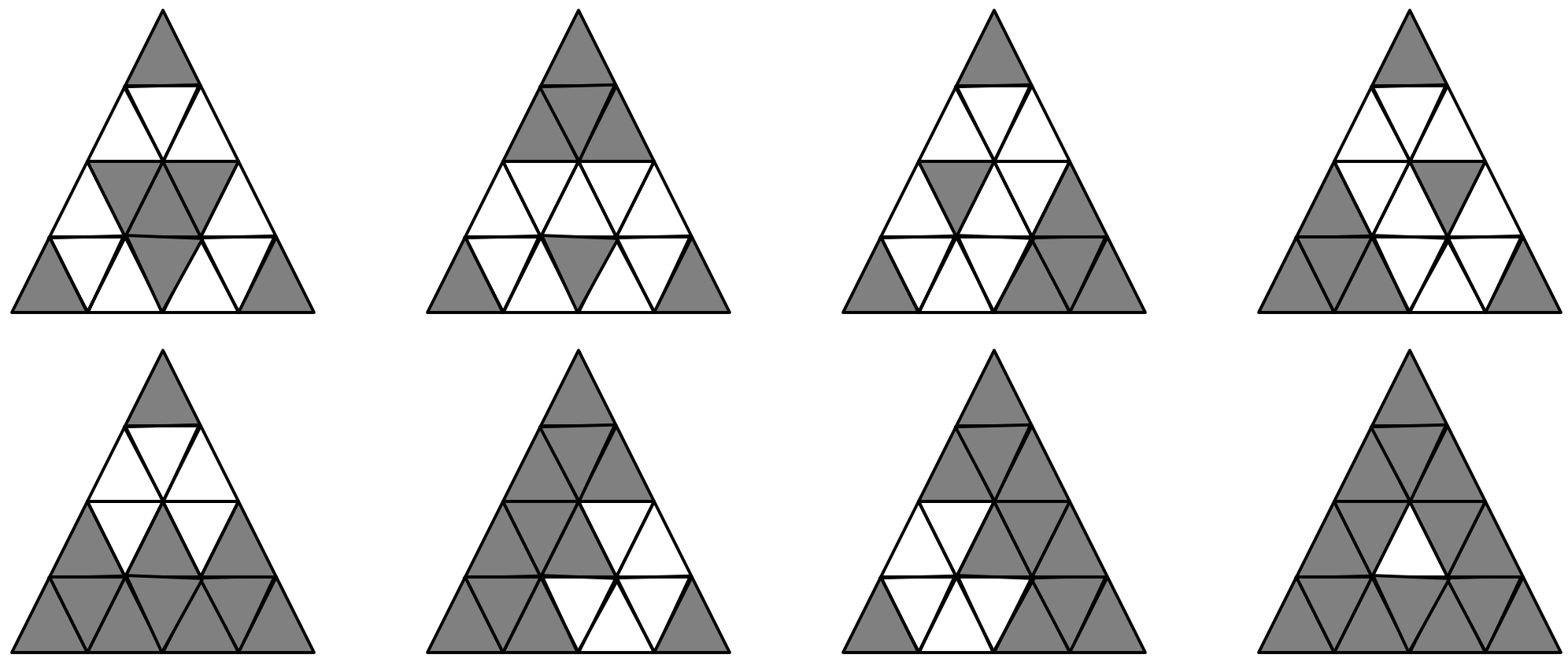

Figure 4.4.  The eight different solutions for the triangular board of five rows.

**Example 4.7.2** Giving the boards a memory value, we can develop Example 4.7.1 further. That is, we can assign a number such that the coins turn from head to tails (or tails to heads) only after a push has acted upon a vertex that number of times; for example 3 times (i.e. when a vertex is labeled 0, 1, or 2 the coin will be heads, and when it is labeled 3, 4, or 5 it will be tails). Then, since $\chi(G) = 3$, we can change one board into the other, as long as, again, $P[L_1(G)] = P[L_2(G)]$.

The boards used for these games can also be structured upon manifolds having a genus other than 0. For example, imagine a 2-simplex graph upon a torus (or Klein bottle), with pushes acting only upon those 2-simplexes whose convex hull (i.e. a 2-simplex in the topological sense of the term) is simply connected. This extends naturally to any compact connected n-manifold. We can therefore imagine (trivially) a 1-simplex graph on a line in 1-dimension; a simple closed curve in 2-dimensions; a knot in 3-dimensions; etc. (And similarly for other n-simplex graphs.)

These games can also be extended to two (or more) person games in a manner such as the following. Given a triangular board made up of heads, as described above (or similarly with a hexagonal board), two players are each assigned a corner. The object of the game is then to form a path of tails to the third corner, with the players alternating pushes and a win accruing for whomever succeeds first.

## 4.8. An Algorithm for (*n*+1)-Colorability

Since if $\chi(G) = n+1$ there are $m^n$ different classes (by Corollary 4.6.4) (given that *G* is a region-connected n-simplex graph), as we move through all $m^r$ different moves possible, we obtain $\frac{1}{m^n}$ of the possible labelings for *G*. If, however, the number of different classes is $m^{n-1}$, then in the worst-case scenario, it is possible that in the first $\frac{1}{m^n}$ of possible moves every labeling of one of the classes is obtained, and in the next $\frac{1}{m^n}$ of possible moves every labeling of a different class is obtained. In which case it would be possible to determine whether or not *G* is (*n*+1)-colorable in at most $\frac{m^r}{m^n} + 1 = m^{r-n} + 1$ moves. (The situation gets only better if the number of different classes for *G* is less than $m^{n-1}$.)

Therefore, we have only to show that if $G$ is not ($n$+1)-colorable then the number of classes of labelings of $G$ decreases by a factor of $m$.

Let $G$ be a graph that is more than ($n$+1)-colorable. Then there is a proper sub-graph of $G$, say $G'$, such that it is maximally ($n$+1)-colorable, with $V(G') = V(G)$. (That is, if we would add any edge from $E(G) - E(G')$ to $E(G')$, $G'$ would no longer be ($n$+1)-colorable.) Therefore, every edge in $E(G) - E(G')$ must connect two vertices of the same color. What results is that the value of the graph (i.e. $P[L(G')]$) would change if acted upon by a push; and the labeling would therefore no longer be of the same class. The actual value would change by $i_k \cdot i_{k+1}^{m-1}$, where $i_k$ is the color of both ends of the new edge. (Without loss of generality we can allow $i_{k+1}$ be the color that is missing.) This is so since a push in this case would result in a change of $i_0 \cdot i_1 \cdots i_k \cdot i_k \cdot i_{k+2} \cdots i_n$. But since $i_0 \cdot i_1 \cdots i_k \cdot i_{k+2} \cdots i_n = i_0 \cdot i_1 \cdots i_k \cdot i_{k+1} \cdot i_{k+1}^{m-1} \cdot i_{k+2} \cdots i_n = (i_0 \cdot i_1 \cdots i_n) i_{k+1}^{m-1} = i_{k+1}^{m-1}$, we have $i_0 \cdot i_1 \cdots i_k \cdot i_k \cdot i_{k+2} \cdots i_n = i_k \cdot i_{k+1}^{m-1}$. However, $i_k \cdot i_{k+1}^{m-1}$ is of order $m$ (where $m$ is the size of the labeling set), so this new edge introduces $m$ new values. Thus the number of classes of labels is decreased by a factor of $m$ (since every class of labelings are now associated to $m$ other classes). For every new edge now added the number of classes could decrease by a factor of $m$, depending upon whether or not each new edge creates a new value when a push is applied to its n-simplex. An ($n$+2)-colorable graph might therefore have $m^{n-1}, m^{n-2}, \ldots, 1$ labeling classes. We have shown that **Lemma 4.8.1** *An (n+2)-colorable (region-connected n-simplex) graph has at most* $m^{n-1}$ *labeling classes.*

For example, the following 4-colorable 2-simplex graph has two labeling classes. 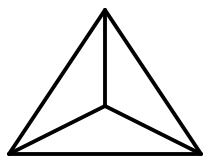

While this 4-colorable 2-simplex graph has only one labeling class. 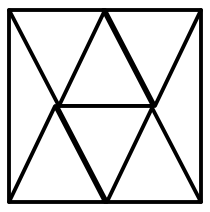
(I.e. it is possible to achieve any labeling of this graph with only pushes.)

Given a graph $G$, we form a new graph, $G'$, with $V(G') = V(G)$ and $E(G') = E(G) + \{$enough edges so that $G'$ would be a region-connected n-simplex graph$\}$. Then

**Corollary 4.8.2** $m^{r-n} + 1$ *moves is sufficient to determine whether a graph G is (n+1)-colorable, where* $r = R(G')$. *If, in addition, G is a region-connected n-simplex graph, then it is necessary as well.*

Since $m$ is independent of $G$ (i.e. dependent upon the labeling set only), it can be chosen arbitrarily. We can therefore improve our bound simply by letting $m = 2$.

Given an n-simplex graph $G = \bigcup_{i=1}^{N} C_i$ with each $C_i$ region-connected and $C_i \bigcap C_j$, $i \neq j$, a subset of an n-simplex, we form a new graph $G'$ by associating to each $C_i \subset G$ a vertex $c_i \subset V(G')$, where $c_i c_j \in E(G')$ given that $C_i \cap C_j$, $i \neq j$, is a subset of an n-simplex.

**Corollary 4.8.3** *Let* $G = \bigcup_{i=1}^{N} C_i$ *be an n-simplex graph where each* $C_i$ *is region-connected and where* $C_i \cap C_j$, $i \neq j$, *is a subset of an n-simplex. Then the condition in Corollary 7.2 is both necessary and sufficient for determining whether or not G is (n+1)-colorable if it's associated graph,* $G'$, *has no cycles.*

**Proof** If $G$ is ($n$+1)-colorable, then given any $C_i$ and $C_j$ as stated, since their intersection is a subset of an n-simplex, the coloring of $C_j$ can be chosen based upon that of $C_i$ (owing to the ($n$+1)-colorability of $C_j$). (The fact that there are no cycles in $G'$ guarantees this ability.) Thus it is possible to choose the coloring of $C_j$ (even if $G$ is planar), such that we can add enough edges into $E(G)$, whose vertices are colored differently, so that $G$ now becomes region-connected without changing the ($n$+1)-colorability of $G$.

While the bounds are given for arbitrary graphs, given specific conditions we can improve upon this bound. For example, for planar graphs, since $r \leq 2v - 4$, we have as a sufficient condition for 3-colorability that we need only try $2^{2v-7} + 1$ moves.

# 5 Open Problems

A number of problems present themselves throughout this paper. These naturally imply different directions we see for further research. Below we have listed some of the main problems both directly related and ancillary to those topics discussed.

1 Is there any direct method of accomplishing knot multiplication (or factoring) directly from a polynomial representation?

2 The secret-key agreement protocol proposed was considered under considerations unique to knot composition and its methods. It is not known whether an attack via alternate representations of knots (as, say, via Dynnikov's three page link diagrams [20] or arc-presentations [19]) would cause a problem.

3 In the arbitrary selection of knots used for the secret-key agreement protocol it is possible that the unknot is used. What complications might arise from this situation?

4 Can second-level distributions (i.e. the distribution of distributions) be quantified in some natural way? And if so, what applications might there be for information theory and cryptology?

5 There lays a relation, hidden though it may be, between games and knots, links, and braids the investigation of which, especially toward the unknotting problem, appears as if it would yield interesting results. A first question toward this end is how might we go about classifying knots via the use of games?

6 How might we use the approaches and techniques already presented to improve upon the upper bounds for $(n+1)$-colorability?

# 6 Related Research Areas

Related to the areas and topics presented we have pursued the interrelation between cryptography, probability, and combinatorics (in particular graph theory) further. We have begun by trying to relate the idea of entropy as introduced by Shannon [55] in the realm of information theory to graphs. Though the idea of entropy as applied to graphs is not new (see Moshowitz [45] and Körner [39]), our approach and definition are different.

Shannon's entropy is given by the following formula: $H(x) = -\sum_{x} p(x)\log_2 p(x)$, where $x$ varies over the entire space.

The idea of Graph Information Entropy we wish to introduce is as follows.

Let $G$ be a graph on $v = V(G)$ vertices and $e = E(G)$ edges. Let the vertices of $G$ be arbitrarily labeled with the values 1 to $n$.

**Definition:** A *Probability Graph* is a graph $G$ with each vertex $i$ assigned a probability $p(i) = \dfrac{\psi(i)}{\sum_{i=1}^{n}\psi(i)}$, where $\psi(i) = \sum_{k=1}^{n}\left(\tfrac{1}{2}\right)^{d(i,k)} - 1$ and $d(x,y)$ is the distance between the vertices $x$ and $y$.

**Definition:** We call $G_\psi = \sum_{i=1}^{n}\psi(i)$ the *$\psi$-value of the graph G.*

**Definition:** The *Graph-Information entropy (GI-entropy value)* of a graph is the value

$$H = -\sum_{i=1}^{n} p(i)\log_2 p(i).$$

We call $\frac{G_\psi}{\binom{n}{2}}$ the *Information-Flow (IF)* of a graph *G*.

We would have, therefore, that a totally unconnected graph would have an IF value of 0, while a complete graph would have a value of 1.